\documentclass[twocolumn]{autart}    
\usepackage{amsfonts}
\usepackage{amsmath}
\usepackage{graphicx}      
\usepackage{float} 
\usepackage{psfrag} 
\usepackage{subfig}
\usepackage{cite}
\usepackage{url}

\begin{document}

\begin{frontmatter}

\title{Sparse solution of the Lyapunov equation for large-scale interconnected systems\thanksref{footnoteinfo}} %

\thanks[footnoteinfo]{Corresponding author: A. Haber, \textit{email address:} aleksandar.haber@gmail.com.}

\vspace{-5mm}
\author{Aleksandar Haber and Michel Verhaegen}    

\address{Delft University of Technology, 2628 CD Delft, The Netherlands}  


\begin{abstract}
We consider the problem of computing an approximate banded solution of the continuous-time Lyapunov equation $\underline{A}\underline{X}+\underline{X}\underline{A}^{T}=\underline{P}$, where the coefficient matrices $\underline{A}$ and $\underline{P}$ are large, symmetric banded matrices. The (sparsity) pattern of $\underline{A}$ describes the interconnection structure of a large-scale interconnected system. Recently, it has been shown that the entries of the solution $\underline{X}$ are spatially localized or decaying away from a banded pattern. We show that the decay of the entries of $\underline{X}$ is faster if the condition number of $\underline{A}$ is smaller. By exploiting the decay of entries of $\underline{X}$, we develop two computationally efficient methods for approximating $\underline{X}$ by a banded matrix. For a well-conditioned and sparse banded $\underline{A}$, the computational and memory complexities of the methods scale linearly with the state dimension. We perform extensive numerical experiments that confirm this, and that demonstrate the effectiveness of the developed methods. The methods proposed in this paper can be generalized to (sparsity) patterns of $\underline{A}$ and $\underline{P}$ that are more general than banded matrices. The results of this paper open the possibility for developing computationally efficient methods for approximating the solution of the large-scale Riccati equation by a sparse matrix.
\end{abstract}
\end{frontmatter}
\section{Introduction}
Large-scale interconnected systems consist of the interconnection of a large number of dynamical subsystems \cite{haber2014Gramian,pakazad2014distributed,dandrea2003,jovbamTAC05platoons,bamjovmitpat12,bamieh02distributedcontrol,khan2008distributing,gorinevsky2008design,motee2008,motee2013measuring,siami2014graph,matni2014communication,andersen2014,zhou2014,zhou2015controllability}. The focus of this paper is on the large-scale interconnected systems described by state-space models with (sparse) banded matrices\footnote{The results of this paper can be extended to matrices with more general sparsity patterns.}. The importance of this class of interconnected systems is best illustrated by the fact that state-space models with banded matrices are obtained by discretizing Partial Differential Equations (PDEs) using the finite difference or finite element methods \cite{benner2004,haberThesis}. Each discretization node or a group of discretization nodes can be seen as a subsystem, and the discretization mesh can be interpreted as a network of dynamical systems\cite{haberThesis,motter2013}. 
\\
The Lyapunov equation is ubiquitous in systems and control theory and in signal processing \cite{gajic2008lyapunov}. For us, the most interesting application of the Lyapunov equation, is in the methods for solving the optimal control (estimation) problems of large-scale systems \cite{mehrmann1991autonomous,simoncini2013computational}. Namely, the solution of the Linear Quadratic (LQ) optimal control problem can be found by solving the Riccati equation. A widely used method for solving the Riccati equation is the Newton method \cite{bini2012numerical,mehrmann1991autonomous,feitzinger2009inexact,wang2014inexact,benner2008numerical,bini2008fast,benner2013numerical}. In each step of the Newton method it is necessary to solve the Lyapunov or Sylvester equations (or the Stein equation for discrete-time systems, see for example \cite{Benner2011}). Further applications of the Lyapunov and Sylvester equations can be found in \cite{simoncini2013computational,gajic2008lyapunov}. 
\\
The goal of this paper is to analyze the solution of the Lyapunov equation for large-scale interconnected systems and to develop computationally efficient methods for approximating the solution by a (sparse) banded matrix. Specifically, we consider the large-scale, continuous-time Lyapunov equation:
\begin{align}
\underline{A}\underline{X}+\underline{X}\underline{A}^{T}=\underline{P}
\label{LyapunovEquation}
\end{align}
where $\underline{A}\in \mathbb{R}^{Nn\times Nn}$ is a (symmetric) negative definite, banded matrix describing the global dynamics of an interconnected system, $\underline{X} \in \mathbb{R}^{Nn\times Nn}$ is a solution that we are searching for, $\underline{P}\in \mathbb{R}^{Nn\times Nn}$ is a banded, negative definite matrix, $N$ is a  large number representing the total number of subsystems and $n\ll N$ is the state order of subsystems. Obviously, finding the solution of \eqref{LyapunovEquation} is a computationally challenging task. The computational challenges in solving the Lyapunov equation for large-scale systems are perhaps best described by the author of the recent survey \cite{simoncini2013computational} on linear matrix equations: \textit{"A distinctive feature in the large-scale setting is that coefficient matrices (the matrices $\underline{A}$ and $\underline{P}$) may be sparse, the solution matrix is usually dense and thus impossible to store in memory... For $\underline{A}$ in the order of $10^4$ or larger the solution cannot be stored explicitly..."}
\par
A large variety of methods for solving the large-scale Lyapunov equation are coping with this problem by searching for a low rank approximation $\tilde{\underline{X}}=ZZ^{T}$ to the "true" solution $\underline{X}$, where $Z$ is the "tall" matrix that is computed and stored \cite{simoncini2013computational}. However, the approximate solution in the form of $\tilde{\underline{X}}=ZZ^{T}$ is a completely dense (fully populated) matrix. Taking into account that the Newton method solves the Riccati equation by solving series of the Lyapunov equations, the solution to the Riccati equation is also dense. This implies that the feedback matrix of the LQ control law is also dense. However, for the distributed control of large-scale interconnected systems, we would like to compute a sparse feedback matrix \cite{motee2013measuring,linfarjovTAC13admm,linfarjovTAC11al}. Namely, a sparse feedback matrix enables us to implement the controller on a network of sensors, actuators and computing units that communicate locally. On the other hand, due to the fact that the computational and memory complexities of multiplying a vector with a sparse matrix are linear, a sparse feedback matrix implies that the centralized LQ control law can be implemented with linear complexity.
\\
If it were possible to accurately approximate the solution of the Lyapunov equation by a sparse matrix, then by using the inexact Newton methods \cite{feitzinger2009inexact} it would be possible to determine a sparse approximate solution of the Riccati equation. That is, if the solution of the Lyapunov equation can be accurately approximated by a sparse matrix, then the LQ feedback matrix can also be approximated by a sparse matrix.  
\\
 In \cite{linfarjovTAC13admm,linfarjovTAC11al,schuler2011}, several methods have been developed for computing sparse, optimal feedback gains. However, the computational and memory complexities of these methods are $O(N^3)$ and $O(N^2)$, respectively, and consequently, these methods are not applicable to large-scale systems.
\\
The above explained problems motivate us to search for the answers to the following questions:
\\\\
\textbf{(1)} Is the solution $\underline{X}$ \textit{spatially localized}? Under the term of a "spatially localized matrix", we understand a matrix whose entries decay quickly in magnitude outside a (sparsity) pattern (for example, outside a bandwidth of a matrix) \cite{benzi2007}. Off-diagonally decaying matrices \cite{benzi2007,Haber:14subspace,demko1984} are typical examples of spatially localized matrices. Roughly speaking, the off-diagonal elements are bounded by an exponential function that decays away from the main diagonal. 
\\
\textbf{(2)} Can the solution of \eqref{LyapunovEquation} be accurately approximated by a banded matrix, possibly sparse\footnote{Under the term of a "sparse banded matrix", we understand a banded matrix whose bandwidth is much smaller than its dimensions.}, and can this approximate solution be computed with $O(N)$ computational and memory complexities? There is a strong correlation between this question and the first one, because if a matrix is spatially localized then it can be accurately approximated by a sparse matrix \cite{haber2014Gramian}.\\\\
Regarding the first question, it has already been shown that $\underline{X}$ is a spatially localized matrix. Namely, the Lyapunov equation can be rewritten as a linear system of equations, in which the coefficient matrix is a Kronecker sum of the matrix $\underline{A}$ and the identity matrix. In \cite{benzi2015decay,canuto2014decay}, it has been shown that the inverses of symmetric, positive definite, banded matrices with a Kronecker sum structure are off-diagonally decaying matrices (with a non-monotonic decay) and several bounds on the off-diagonal decay rate have been derived. These important results imply that $\underline{X}$ is a spatially localized matrix \cite{simoncini2015lyapunov}. However, because the upper-bounds on the off-diagonal decay rate presented in \cite{benzi2015decay,canuto2014decay,simoncini2015lyapunov} have a relatively complex integral form, they need to be evaluated numerically. 
Consequently, from these bounds it might be hard to draw important conclusions on how the condition number of banded $\underline{A}$ influences the decay rate of entries of $\underline{X}$. It should come as no surprise that the condition number of $\underline{A}$ dominantly determines the decay rate of $\underline{X}$, because it is well-known that the decay rates of matrix functions, such as $\underline{A}^{-1}$, are primarily determined by the condition number of $\underline{A}$ \cite{demko1984}. On the other hand, in \cite{benzi2007,benzi2015decay}, it has been shown that the matrix exponential of $\underline{A}$, denoted by $\text{exp}\left(t\underline{A}\right)$, where $t$ is time, is also an off-diagonally decaying matrix. 
\par
By exploiting the fact that the solution of the Lyapunov equation $\underline{X}$ and the matrix exponential $\text{exp}\left(\underline{tA}\right)$ are spatially localized matrices, we develop two computationally efficient methods for approximating $\underline{X}$ by a banded matrix. Furthermore, we show that the decay of entries of $\underline{X}$ is faster if the condition number of $\underline{A}$ is smaller. Our results indicate that for a well-conditioned and sparse banded $\underline{A}$, the proposed methods are able to approximate the matrix $\underline{X}$ with $O(N)$ complexity. The results of this paper open the possibility for developing computationally efficient methods for approximating the solution of the large-scale Riccati equation by a sparse matrix. Furthermore, the results of this paper can also be generalized to diagonalizable banded matrices $\underline{A}$, and to matrices with more general (sparsity) patterns, see Section \ref{generalSparsityPattern}.
\par
This paper is organized as follows. In Section \ref{problemFormulation}, we present the problem formulation. In Section \ref{newDecaySection}, we analyze the decay of entries of $\underline{X}$ and we present a method for estimating the a priori pattern of an approximate solution. In Section \ref{methodsApproximation}, we develop two methods for approximating $\underline{X}$. In Section \ref{sectionNumericalExperiments}, we present numerical experiments, and in Section \ref{sectionConclusions} we present conclusions.
\section{Problem formulation}
\label{problemFormulation}
\subsection{Notation}
\label{notationSection}
The notation $X=[x_{i,j}]$ denotes a matrix whose $(i,j)$ entry is $x_{i,j}$, whereas  $X=[X_{i,j}]$ denotes a block matrix whose $(i,j)$ entry is the matrix $X_{i,j}$. 
The notation $\mathbf{z}=\text{col}\left(\mathbf{z}_{1},\mathbf{z}_{2},\ldots,\mathbf{z}_{M} \right)$ stands for 
$\mathbf{z}=[\mathbf{z}_{1}^{T} \;\mathbf{z}_{2}^{T}\; \ldots \; \mathbf{z}_{M}^{T} ]^{T}$. An $N\times N$ matrix $X=[x_{i,j}]$ is called a \textit{banded matrix} if there exists an even positive integer $s$, such that $x_{i,j}=0$ when $|i-j|>s/2$ \cite{demko1984}. The number $s$ is called the \textit{bandwidth} of $X$ and we say that the matrix $X$ is \textit{s-banded} \cite{demko1984}. For example, a tridiagonal matrix is $2$-banded. If $s\ll N$, then $X$ is called a \textit{sparse banded matrix}. If the matrix $X_{1}$ has the bandwidth $s_{1}$ and the matrix $X_{2}$ has the bandwidth $s_{2}$, then the product $X_{1}X_{2}$ has the bandwidth equal to $s_{1}+s_{2}$.
The notations $\left\|X \right\|_{2}$ and $\left\| X\right\|_{F}$ denote the 2-norm and the Frobenius norm of $X$, respectively. The symbol $\otimes$ denotes the Kronecker product, and the operator $\text{vec}\left(X \right)$ is a standard "vec" operator \cite{verhaegen2007}. The (column) vector $\mathbf{q}_{j}$ denotes a vector that has all zeros except $1$ on the position $j$ \cite{simoncini2015lyapunov}. For example, a matrix $Z$, having only one non-zero element $z_{i,j}$ at position $(i,j)$, can be represented by $Z=\mathbf{q}_{i}z_{i,j}\mathbf{q}_{j}^{T}$.
\par
We consider a subsystem $\mathcal{S}_{i}$:
\begin{align}
\mathcal{S}_{i} \left\{ \begin{array} {rll}
\dot{\mathbf{w}}_{i}(t)&=A_{i,i}\mathbf{w}_{i}(t)+\sum_{j=i-b, j\ne i }^{i+b}A_{i,j}\mathbf{w}_{j}(t) \\ 
\mathbf{y}_{i}(t)&=C_{i}\mathbf{w}_{i}(t)   \end{array} \right.
\label{localSubSys}
\end{align}
where $\mathbf{w}_{i}(t)\in \mathbb{R}^{n}$ is the \textit{local state} of the subsystem $\mathcal{S}_{i}$ and $\mathbf{y}_{i}(t)\in \mathbb{R}^{r}$ is the \textit{local output}, $A_{i,j} \in \mathbb{R}^{n\times n}$ and $C_{i} \in \mathbb{R}^{r\times n}$ and $t$ is time. The state-space model of the \textit{global system} $\mathcal{S}$ is:
\begin{align}
\mathcal{S} \left\{ \begin{array} {rl}
\dot{\underline{\mathbf{w}}}(t)&=\underline{A}\underline{\mathbf{w}}(t)\\
\underline{\mathbf{y}}(t)&=\underline{C}\underline{\mathbf{w}}(t) \end{array} \right. 
\label{globalSys}
\end{align}
where $\underline{\mathbf{w}}(t)=\text{col}\left(\mathbf{w}_{1}(t), \ldots, \mathbf{w}_{N}(t)\right)$, \\ $\underline{\mathbf{y}}(t)=\text{col}\left( \mathbf{y}_{1}(t),\ldots,\mathbf{y}_{N}(t)\right)$, $\underline{A}\in \mathbb{R}^{Nn\times Nn}$ and $\underline{C}\in \mathbb{R}^{Nr\times Nn}$. The vectors $\underline{\mathbf{w}}(t)\in \mathbb{R}^{Nn}$ and $\underline{\mathbf{y}}(t)\in \mathbb{R}^{Nr}$ are called the \textit{global state} and \textit{global output}, respectively. We assume that the total number of subsystems $N$ is a large number and that $n\ll N$. Furthermore, we assume that $b\ll N$. That is, we assume that the matrix $\underline{A}$ is a sparse banded matrix. The bandwidth of $\underline{A}$ is denoted by $m\ll N$. We assume that the matrix $\underline{A}$ is symmetric and asymptotically stable (although all the methods in this paper can be generalized to banded diagonalizable matrices and to matrices with more general patterns, see Section \ref{generalSparsityPattern}). Finally, we assume that the matrix $\underline{P}$ in \eqref{LyapunovEquation} is a sparse, banded matrix with the bandwidth equal to $l$, where $l\ll N$. For example, the model \eqref{globalSys} can be obtained by discretizing the 2D or 3D partial differential equations using the finite difference methods \cite{benner2004,haberThesis}.
\par
Methods for solving \eqref{LyapunovEquation} will be built upon on the following two representations of the solution $\underline{X}$ \cite{simoncini2015lyapunov,laub2005matrix}.  Taking into an account the symmetry of $\underline{A}$, the (unique) solution of  \eqref{LyapunovEquation} has the integral representation \cite{laub2005matrix}: 
\begin{align}
\underline{X}=-\int_{0}^{\infty}\exp\left(t\underline{A}\right)\underline{P} \exp\left(t\underline{A}\right) \mathrm{d}t
\label{integralRepresentationLyapunov}
\end{align}
 By vectorizing \eqref{LyapunovEquation}, we obtain: 
\begin{align}
\mathcal{A}\overline{\mathbf{x}}&=\overline{\mathbf{p}}, \label{KroneckerRepresentation} \\ \overline{\mathbf{x}}&=\mathcal{A}^{-1}\overline{\mathbf{p}} \label{KroneckerRepresentation2}
\end{align}
where $\mathcal{A}\in \mathbb{R}^{ (Nn)^{2}\times (Nn)^{2}}$, $\overline{\mathbf{x}},\overline{\mathbf{p}} \in \mathbb{R}^{(Nn)^2}$ are defined by:
\begin{align}
\mathcal{A}=I \otimes \underline{A}+\underline{A}\otimes I,\; \overline{\mathbf{x}}=\text{vec}\left(\underline{X} \right), \; \overline{\mathbf{p}}=\text{vec}\left(\underline{P} \right)
\label{KroneckerRepresentation3}
\end{align}
and where $I$ is an $Nn\times Nn$ identity matrix. Based on the fact that the entries of the solution $\underline{X}$ are decaying away from a banded pattern (in some cases they exhibit a damped oscillatory behavior) \cite{benzi2015decay,canuto2014decay,simoncini2015lyapunov}, in the sequel we will develop two methods for approximating $\underline{X}$ by a (sparse) banded matrix. The first method looks for an approximate solution by solving a least-squares problem formed on the basis of \eqref{KroneckerRepresentation}. The least-squares problem is formed by eliminating the columns of $\mathcal{A}$ and the elements of $\overline{\mathbf{x}}$ corresponding to the small elements of $\underline{X}$ that are predicted by the decaying property. The second method consists of the following two steps. In the first step, the integral in \eqref{integralRepresentationLyapunov} is approximated by a banded matrix. This approximation is obtained by exploiting the fact that $\exp\left(t\underline{A}\right)$ is an off-diagonally decaying matrix \cite{benzi2007,benzi2015decay}. In the second step, the accuracy of this approximate solution is additionally improved by using the gradient projection method \cite{bertsekas1999nonlinear}. We start with the analysis of the decaying behavior of entries of $\underline{X}$. 
\section{Decay rate analysis and a priori pattern of $\underline{X}$}
\label{newDecaySection}
The goal of this section is to analyze the influence of the condition number of $\underline{A}$ on the decay of entries of $\underline{X}$, and to develop a method for the prediction of the \textit{a priori (sparsity) pattern} of the approximate solution. Important insights obtained in this section will be used in Section \ref{methodsApproximation} to develop computationally efficient approximation methods. 
\begin{defn}
\label{offDiagonalDecayDefinition} \cite{benzi2007,demko1984,canuto2014decay,shao2014finite}
We say that an $Nn\times Nn$ matrix $Z=[z_{i,j}]$ is an off-diagonally decaying matrix if there exist $\tau >0$ and $\rho\in (0,1)$ such that $|z_{i,j}| \le \tau \rho^{|i-j|}$
for all $i,j=1,\ldots,Nn$. $\hfill \square$ \end{defn}
\par
The constant $\rho$ is referred to as the \textit{decay rate} of $Z$ \cite{shao2014finite}.  We will use theoretical results developed in \cite{demko1984,benzi1999bounds,benzi2007,canuto2014decay,benzi2015decay} to analyze the decay rate of $\underline{X}$, as well as to provide insights on how the condition number of $\underline{A}$ influences the decay rate.
\subsection{The influence of the condition number of $\underline{A}$ on decay of $\underline{X}$}
\label{sectionDecayingX}
The constants $a$ and $b$ are defined by $a=\lambda_{\text{min}}\left(\underline{A} \right)$ and $b=\lambda_{\text{max}}\left(\underline{A} \right)$, where $\lambda_{\text{min}}
\left(\cdot\right)$ and $\lambda_{\text{max}}\left(\cdot \right)$ denote minimal and maximal eigenvalues, respectively. Given that the matrix $\underline{A}$ is a symmetric, asymptotically stable matrix, we have that $a<0$ and $b<0$ and $|b|\le |a|$. The complexity analysis of computing $a$ and $b$ is discussed in Section \ref{complexityAnalysis}. The condition number of $\underline{A}$ will be denoted by $\kappa$. It follows that $\kappa=a/b$.
\par
For presentation clarity, we will first consider a special case, when the matrix $\underline{P}$ is a diagonal matrix, and subsequently we will consider a more general case when the matrix $\underline{P}$ is a banded matrix (or even fully populated). Let us suppose that $\underline{P}=\gamma I$, $\gamma <0$. Then, because $\underline{A}$ is symmetric, it can be verified that $\underline{X}=\left((2/ \gamma ) \underline{A} \right)^{-1}$ satisfies \eqref{LyapunovEquation}. Due to the fact that $\gamma <0$, we have $ \underline{X}=\left( (2/\gamma) \underline{A} \right)^{-1}=(|\gamma|/2) \left(-\underline{A} \right)^{-1} $. Because by assumption $\underline{A}$ is asymptotically stable, symmetric matrix, the matrix $-\underline{A}$ is positive definite, and consequently, from Theorem 2.4 in \cite{demko1984}, it follows that $\left(-\underline{A} \right)^{-1}$ is an off-diagonally decaying matrix. This further implies that $X=[x_{i,j}]$ is an off-diagonally decaying matrix, that is,  $|x_{i,j}|\le \tau \rho^{|i-j|}$, where 
\begin{align}
&\tau =\frac{|\gamma|}{2} K_{1},\;\;\; K_{1}=\frac{1}{|b|} \max{\{1,\frac{\left(1+\sqrt{\kappa} \right)^{2}}{2 \kappa} \}}, \notag \\ & \rho=\left( \frac{\sqrt{\kappa}-1}{\sqrt{\kappa}+1}\right)^{\frac{2}{m}}
\label{decayRateLyapunov}
\end{align}
From \eqref{decayRateLyapunov} we see that the decay rate of $\underline{X}$ depends on the condition number of $\underline{A}$. Specifically, if $\underline{A}$ is well-conditioned ($\kappa$ is close to 1), then the decay rate\footnote{Beside $\kappa$, we see that the minimal singular value $|b|$ of $\underline{A}$ determines the off-diagonally decaying behavior. Throughout the paper we will assume that $|b|$ is not very small, which implies that the decaying behavior of $\underline{A}$ is primarily determined by $\kappa$.} is fast (the number $\rho$ is small). Now, does $\underline{X}$ exhibit a similar behavior when $\underline{P}$ is a sparse banded matrix? Not surprisingly, the answer is yes. Before we show this, it should be first observed that the matrix $\mathcal{A}$ is also an off-diagonally decaying matrix. Namely,  from \cite{laub2005matrix} (Theorem 13.16) it follows that:
\begin{align}
\lambda_{\text{max}}\left(\mathcal{A}\right)=2\lambda_{\text{max}}\left(\underline{A}\right)=2b,\;
\lambda_{\text{min}}\left(\mathcal{A}\right)=2\lambda_{\text{min}}\left(\underline{A}\right)=2a
\label{tmp111}
\end{align}
This implies that the condition numbers of $\mathcal{A}$ and $\underline{A}$ are equal. On the other hand, because $\underline{A}$ is $m$-banded, the matrix $\mathcal{A}$ is $m_{1}$-banded, with $m_{1}=Nnm$ \cite{canuto2014decay}. Furthermore, the matrix $\mathcal{A}$ is negative definite.  Similarly to the analysis of the decay rate of $\underline{X}$ (for the case of a diagonal $\underline{P}$), applying the results of Theorem 2.4 in \cite{demko1984} to $-\mathcal{A}$, we conclude that $\mathcal{A}^{-1}$ is an off-diagonally decaying matrix, with the decay rate specified by:
\begin{align}
&\tau_{1} =\frac{1}{2|b|} \max{\{1,\frac{\left(1+\sqrt{\kappa} \right)^{2}}{2 \kappa} \}}, \; \rho_{1}=\left( \frac{\sqrt{\kappa}-1}{\sqrt{\kappa}+1}\right)^{\frac{2}{m_{1}}}
\label{decayRateLyapunov2}
\end{align}
By comparing \eqref{decayRateLyapunov} and \eqref{decayRateLyapunov2}, we conclude that the decay rates $\rho$ of $\underline{X}$ (for the case of diagonal $\underline{P}$) and $\rho_{1}$ of $\mathcal{A}^{-1}$, have the same exponential base (determined by $\kappa$). Furthermore, we see that because $m_{1}$ is by several orders of magnitude larger than $m$, the decay rate of $\mathcal{A}^{-1}$ is slower than the decay rate of $\underline{X}$. However, the matrix $\mathcal{A}^{-1}$ is much larger than $\underline{X}$. Furthermore, taking into account the sizes and the bandwidths of the corresponding matrices, the entries of both matrices that are far away from the corresponding main diagonals or bandwidths, are small when the matrix $\underline{A}$ is well-conditioned. To further analyze the decay rate of $\underline{X}$, we represent the matrix $\underline{P}=[p_{i,j}]$ as follows \cite{simoncini2015lyapunov}:
\begin{align}
\underline{P}=\sum_{i=1}^{Nn}\sum_{j= 1}^{Nn} \underline{P}_{i,j}
\label{expansionP}
\end{align}
where $\underline{P}_{i,j}=\mathbf{q}_{i}p_{i,j}\mathbf{q}_{j}^{T}$ and where $\mathbf{q}_{i}$ and $\mathbf{q}_{j}$ are defined in Section \ref{notationSection}. For notation simplicity and presentation clarity, in \eqref{expansionP} and throughout the rest of this section, we have formally ignored the fact that most of the entries of $\underline{P}$ are zero (the matrix $\underline{P}$ is $l$-banded, so its entries $p_{i,j}$ for which $|i-j|>l/2$, are equal to zero). While interpreting the results, it should be kept in mind that the summation in \eqref{expansionP} and in subsequent expressions, should be performed only with respect to the indices $(i,j)$ belonging to the bandwidth region. Let $\underline{X}_{i,j}\in \mathbb{R}^{Nn\times Nn}$ be the solution of the Lyapunov equation:
\begin{align}
\underline{A}\underline{X}_{i,j}+\underline{X}_{i,j}\underline{A}=\underline{P}_{i,j}=\mathbf{q}_{i}p_{i,j}\mathbf{q}_{j}^{T}
\label{decomposedLyapunov}
\end{align}
where $i,j=1,\ldots,Nn$. Then, because of the linearity of the Lyapunov equation it follows that its solution can be decomposed as follows \cite{simoncini2015lyapunov}:
\begin{align}
\underline{X}=\sum_{i=1}^{Nn}\sum_{j=1}^{Nn} \underline{X}_{i,j}
\label{LyapunovEquationDecomposition}
\end{align}
where each of $ \underline{X}_{i,j}$ is the solution of \eqref{decomposedLyapunov}. The decomposition \eqref{LyapunovEquationDecomposition} enables us to prove the following theorem.
\begin{thm}
Let $\overline{\mathbf{x}}_{i,j}=\text{vec}\left(\underline{X}_{i,j} \right)$ and let the element of $\overline{\mathbf{x}}_{i,j}$ on the position $s$, $s=1,\ldots,(Nn)^{2}$, be denoted by $\overline{x}_{s}^{i,j}$, then
\begin{align}
|\overline{x}_{s}^{i,j}|\le |p_{i,j}| \tau_{1}\rho_{1}^{|\phi(i,j)-s|}
\label{proofDecay112}
\end{align}
where $\phi(i,j)=(j-1)Nn+i$ and $\rho_{1}$ and $\tau_{1}$ are defined in \eqref{decayRateLyapunov2}. Moreover, let the element of  $\overline{\mathbf{x}}=\text{vec}\left(\underline{X}\right)$ on the position $s$ be denoted by $\overline{x}_{s}$, then
\begin{align}
|\overline{x}_{s}| \le \tau_{1} \sum_{i=1}^{Nn}\sum_{j=1}^{Nn}|p_{i,j}| \rho_{1}^{|\phi(i,j)-s|}
\label{constantDerivation}
\end{align}
\label{theoremLastXdecayProof}
\end{thm}
\vspace{-4mm}
\textit{\textbf{Proof}:} After applying $\text{vec}\left(\cdot\right)$ operator to \eqref{decomposedLyapunov}, we obtain:
\begin{align}
\mathcal{A}\overline{\mathbf{x}}_{i,j}=p_{i,j}\mathbf{q}_{(j-1)Nn+i},\;\;
\overline{\mathbf{x}}_{i,j}=p_{i,j}\mathcal{A}^{-1}\mathbf{q}_{(j-1)Nn+i}
\label{KroneckerRepresentationDecomposed}
\end{align}
 and according to our notation the vector $\mathbf{q}_{(j-1)Nn+i}$ has $1$ on the position $(j-1)Nn+i$. Let the $\phi(i,j)=\left((j-1)Nn+i\right)$-th column of $\mathcal{A}^{-1}$ be denoted by $\boldsymbol{\theta}_{(i,j)} \in \mathbb{R}^{(Nn)^{2}}$. Let the element of $\boldsymbol{\theta}_{(i,j)}$ on the position $s$ be denoted by $\theta^{(i,j)}_{s}$. From \eqref{KroneckerRepresentationDecomposed} it follows that $\overline{\mathbf{x}}_{i,j}$ is equal to the $\phi(i,j)$-th column of $\mathcal{A}^{-1}$ multiplied by a constant $p_{i,j}$:
\begin{align}
\overline{\mathbf{x}}_{i,j}=p_{i,j}\boldsymbol{\theta}_{(i,j)}
\label{proofDecay1112}
\end{align}
Given that $\mathcal{A}^{-1}$ is an off-diagonally decaying matrix, the absolute values of elements of its each column are bounded by an exponential function that decays away from the element on the main diagonal. The element of $\boldsymbol{\theta}_{(i,j)}$ that is on the main diagonal of $\mathcal{A}^{-1}$ is the element $\theta^{(i,j)}_{r}$ for which $r=\phi(i,j)$. All this implies that for the entries on the $\phi(i,j)$-th column of $\mathcal{A}^{-1}$ we can write:
\begin{align}
|\theta_{s}^{(i,j)}|\le \tau_{1}\rho_{1}^{|\phi(i,j)-s|}
\label{proofDecay111}
\end{align}
 From \eqref{proofDecay1112} and \eqref{proofDecay111} we obtain \eqref{proofDecay112}. By vectorizing \eqref{LyapunovEquationDecomposition} and using \eqref{proofDecay112}, we can similarly prove \eqref{constantDerivation}. $\hfill \square$
\par
Although conservative compared to the bounds in \cite{benzi2015decay,canuto2014decay,simoncini2015lyapunov}, the bounds in  \eqref{proofDecay112} and \eqref{constantDerivation} can be used to analyze the dependence of the decay rate of $\underline{X}$ on the condition number of $\underline{A}$. Namely, the upper bound \eqref{proofDecay112} on the entries $\overline{x}_{s}^{i,j}$ of the vector $\overline{\mathbf{x}}_{i,j}$, has a maximum value for the entry $s=\phi(i,j)$. When, on the other hand, $\overline{\mathbf{x}}_{i,j}$ is transformed back to a matrix format, this maximum corresponds to the $(i,j)$ entry of $\underline{X}_{i,j}$. Roughly speaking, we also see that the entries of $\underline{X}_{i,j}$ that are further away from the entry $(i,j)$, are bounded by a function that decays as the distance between these entries and the entry $(i,j)$ is increased. The decay rate is determined by the condition number of $\mathcal{A}$ (or equivalently by the condition number of $\underline{A}$). If $\underline{A}$ is well-conditioned, this decay is fast. Now, taking into account that $\underline{P}$ is sparse and banded, we have that the upper bounds have maximal values for the entries $(i,j)$ located in the bandwidth region of $\underline{P}$. On the other hand, from \eqref{LyapunovEquationDecomposition}, we see that the decaying behavior of $\underline{X}$ is determined by the sum of locally decaying behaviors of each individual $\underline{X}_{i,j}$ (it should be remembered that the sum in \eqref{LyapunovEquationDecomposition} is only performed over the indices $(i,j)$ for which $p_{i,j}$ is not zero, that is, for the entries inside the bandwidth of $\underline{P}$). This leads us to the conclusion that the elements of $\underline{X}$ that are far away from the bandwidth of $\underline{P}$ should be relatively small for well conditioned $\underline{A}$. This insights and the insights provided in Section \ref{generalSparsityPattern}, will enable us to develop computationally efficient algorithms for approximating $\underline{X}$ in Section \ref{methodsApproximation}.
\subsection{Predicting the sparsity pattern of $\underline{X}$}
\label{generalSparsityPattern}
In Section \ref{methodsApproximation}, we will develop computationally efficient methods for approximating $\underline{X}$. To develop these methods, we need to chose an \textit{a priori (sparsity) pattern} of $\underline{X}$. From the previous discussion, we may conclude that if the matrix $\underline{A}$ is well-conditioned then the entries of  $\underline{X}$ that are far away from the bandwidth region of $\underline{P}$ are small. Consequently, the a priori pattern can be chosen as a (sparse) banded matrix. By neglecting the entries of $\overline{\mathbf{x}}$ that are outside this pattern, we can reduce the dimension of the linear system \eqref{KroneckerRepresentation}, and we can obtain its solution by solving a least squares problem. Furthermore, the reduced coefficient matrix of this system is sparse, and this sparsity can be exploited to quickly compute the solution.
\par
However, a banded a priori pattern can sometimes be suboptimal in the sense that the "true" solution $\underline{X}$ has a large number of small entries inside of the banded a priori pattern (this can happen for example when the entries of $\underline{X}$ exhibit an oscillatory behavior). This can lead to unnecessary increase of the computational and memory complexities of the proposed algorithms. Is there a more optimal a priori pattern, that can more accurately capture the "true" behavior of the entries of $\underline{X}$? Moreover, we might ask ourselves what should be an a priori pattern in the case of more general patterns of the coefficient matrices $\underline{A}$ and $\underline{P}$?
\par
Due to the Kronecker sum structure of $\mathcal{A}$, the entries of $\mathcal{A}^{-1}$ exhibit a dominantly oscillatory behavior \cite{simoncini2015lyapunov,canuto2014decay,benzi2015decay}. Relatively good estimates of this behavior are presented in \cite{simoncini2015lyapunov,canuto2014decay,benzi2015decay}. Moreover, from \cite{simoncini2015lyapunov,canuto2014decay,benzi2015decay} it follows that for problems for which $\underline{A}$ is a Kronecker sum, $\underline{A}=I\otimes Z_{1}+ Z_{1} \otimes I$, where $ Z_{1}$ is a sparse banded matrix, the entries of $\underline{X}$ also exhibit an oscillatory behavior. These important insights can be used to design a priori patterns that optimally (meaning that the number of non-zero elements is smaller than the number of non-zero elements of a banded matrix) capture the "true" behavior of the entries of $\underline{X}$. For example, these results can indicate that the a priori pattern of $\underline{X}$ should be a multi-banded matrix (matrix that has series of zero and non zero diagonals below or above the main diagonal, examples are shown in Fig. \ref{fig:SurfSparsity3Dheat}(b) and Fig. \ref{fig:AccErr3Dheat}(d)). On the other hand, decay bounds available in the literature indicate that the results of this paper can be generalized to non-symmetric matrices $\underline{A}$ and $\mathcal{A}$. Namely, Theorem 3.5 in \cite{benzi2007} proves that functions of diagonalizable banded matrices also exhibit a form of the off-diagonal decay. Moreover, from Theorem 3.4 in \cite{benzi2007}, it follows that functions of more general class of diagonalizable sparse matrices (not necessarily banded) are spatially localized.
The practical potential of the results presented in \cite{simoncini2015lyapunov,canuto2014decay,benzi2015decay} for determining the optimal a priori pattern of $\underline{X}$ will be investigated in our future work. 
\par
In this paper, we will present a relatively simple approach that can give us additional insights into a more optimal a priori pattern of $\underline{X}$. Furthermore, this approach can be applied to a more general class of sparse matrices $\underline{A}$ and $\underline{P}$, such as sparse multi-banded matrices or even non-symmetric matrices. Namely, using the Neumann representation of the matrix inverse or the fact that a matrix satisfies its characteristic polynomial, in \cite{huckle1999apriori} it has been shown that a relatively good guess of the a priori pattern of $\mathcal{A}^{-1}$ is given by the pattern of the following matrix:
\begin{align}
\mathcal{C}=I+\mathcal{A}+\mathcal{A}^{2}+\ldots+\mathcal{A}^{z_{1}}
\label{aprioriSeriesExpansion}
\end{align}
For simplicity and without the loss of generality, in \eqref{aprioriSeriesExpansion} we will assume that zero entries of $\mathcal{C}$ are not created by incidental cancellations of non-zero entries of matrices produced by taking and summing up powers of $\mathcal{A}$ \cite{huckle1999apriori}. If $\mathcal{A}$ is $m_{1}$-banded, then from \eqref{aprioriSeriesExpansion} it follows that $\mathcal{C}$ is $z_{1}m_{1}$-banded. Next, if the matrix $\mathcal{A}$ is (sparse) banded, then the matrix $\mathcal{C}$ will also be banded, for small values of $z_{1}$. Also, under the same conditions, multi-banded structure of $\underline{A}$ or consequently of $\mathcal{A}$, will be preserved in $\mathcal{C}$. The fundamental question that needs to be asked is how large $z_{1}$ should be such that the majority of the significant entries of $\mathcal{A}^{-1}$ are captured by the sparsity pattern of $\mathcal{C}$, and is there a connection between $z_{1}$ and the condition number of $\mathcal{A}$? The answers to these questions can be obtained by analyzing the accuracy of the expansions from which the expression \eqref{aprioriSeriesExpansion} originates. However, this might be a nontrivial problem. Instead, we use an alternative way to show that for well-conditioned $\underline{A}$, a relatively good guess of the a priori pattern of $\mathcal{A}^{-1}$ is formed by summing up a relatively low powers of $\mathcal{A}$ in \eqref{aprioriSeriesExpansion}. This will be shown by considering the Newton-Schultz iteration \cite{pan1991improved} for approximating\footnote{The Newton-Schultz iteration will only be used to argue about the relationship between the condition number and the a priori pattern, and it will not be used to compute the a priori pattern.} $\mathcal{A}^{-1}$. The Newton-Schultz iteration is defined by \cite{pan1991improved}:
\begin{align}
\mathcal{B}_{k+1}=\mathcal{B}_{k}\left(2I-\mathcal{A}\mathcal{B}_{k}\right), \; k=0,1,2,\ldots
\label{newtonSchultz}
\end{align}
where $\mathcal{B}_{k}$ is the approximation of $\mathcal{A}^{-1}$ at the $k$-th iteration. The iteration is initialized by $\mathcal{B}_{0}$:
\begin{align}
\mathcal{B}_{0}=\frac{2}{a^2_{1}+b^2_{1}}\mathcal{A}^{T}
\label{initialNewton}
\end{align}
where $a_{1}$ and $b_{1}$ are the maximal and minimal singular values of $\mathcal{A}$, respectively. The accuracy at the $k$-th iteration of the Newton-Schultz iteration is measured by the norm of $\mathcal{E}_{k}=I-\mathcal{A}\mathcal{B}_{k}$. Using the fact that $a_{1}=2|a|$ and $b_{1}=2|b|$ (see \eqref{tmp111}), and using the results of  \cite{pan1991improved}, it can be easily shown that (see for example Theorem 3.1 in \cite{haber2014Gramian})
\begin{align}
\left\|\mathcal{E}_{k}\right\|_{2}\le \left(\frac{\kappa^2-1}{\kappa^2+1}\right)^{2^{k}}
\label{errorFinalBound}
\end{align}
Similarly to \eqref{aprioriSeriesExpansion}, the Newton-Schultz iteration \eqref{newtonSchultz} tells us that the relatively good guess of the a priori pattern of $\mathcal{A}^{-1}$ is given by the sum of powers of $\mathcal{A}$. This can be easily shown by starting from the initial guess \eqref{initialNewton}, and by propagating the recursion \eqref{newtonSchultz}. It should be remembered that from the (sparsity) pattern point of view, in which we only focus on the structure but not on the exact numerical values of the entries, the patterns produced by \eqref{aprioriSeriesExpansion} and \eqref{newtonSchultz} can be made to be identical, and are mainly determined by the maximal powers of $\mathcal{A}$ that are determined by the parameters $z_{1}$ and $k$. Furthermore, the bound \eqref{errorFinalBound} tells us that if the matrix $\underline{A}$ is well-conditioned (it should be remembered that the condition numbers of $\underline{A}$ and $\mathcal{A}$ are equal), then the number $k$ that produces a good approximation accuracy of the Newton-Schultz iteration is relatively small. Under the condition that $\mathcal{A}$ is sparse multi-banded, small $k$ or equivalently, small $z_{1}$, produces sparse multi-banded $\mathcal{C}$. That is, the parameter $z_{1}$ in \eqref{aprioriSeriesExpansion}, that primarily determines the structure of the a priori pattern, should be small for well-conditioned $\underline{A}$. 
\par
Taking into account \eqref{KroneckerRepresentation2} and \eqref{aprioriSeriesExpansion}, we have that the pattern of the vectorized solution $\overline{\mathbf{x}}$ can be estimated by the pattern (or non-zero entries) of the vector $\overline{\mathbf{x}}_{2}$:
\begin{align}
\overline{\mathbf{x}}_{2}&=\mathcal{C}\overline{\mathbf{p}} \label{KroneckerRepresentation22}
\end{align}
For a large $N$, the vector $\overline{\mathbf{x}}_{2}$ cannot be computed directly from 
\eqref{KroneckerRepresentation22}, simply because it is impossible to explicitly form the powers of $\mathcal{A}$. Instead, the vector $\overline{\mathbf{x}}_{2}$ should be expressed as the sum of $\mathcal{A}^{l}\overline{\mathbf{p}}$, $l=1,2,
\ldots, z_{1}$. Then, using the main properties of the Kronecker sum, these terms should be written as matrices (by reversing the vectorizing operation that is used to form \eqref{KroneckerRepresentation2}), and summed up together to form the matrix representation of $\overline{\mathbf{x}}_{2}$. For example, the term $\mathcal{A}\overline{\mathbf{p}}$, when transformed back to a matrix format, has the following form $\underline{A}\underline{P}+\underline{P}\underline{A}$. Similarly, the term $\overline{\mathbf{s}}_{2}=\mathcal{A}^{2}\overline{\mathbf{p}}$ can be written as
\begin{align}
\overline{\mathbf{s}}_{2}=(I \otimes \underline{A}+\underline{A}\otimes I)\overline{\mathbf{p}}_{1}, \;\;\; \overline{\mathbf{p}}_{1}= (I \otimes \underline{A}+\underline{A}\otimes I)\overline{\mathbf{p}}
\label{kroneckerMatrix1}
\end{align}
and when transformed back to a matrix format $\underline{S}_{2}$, the term $\overline{\mathbf{s}}_{2}$ has the following form:
\begin{align}
\underline{S}_{2}=\underline{A}\underline{P}_{1}+\underline{P}_{1}\underline{A},\;\;\; \underline{P}_{1}=\underline{A}\underline{P}+\underline{P}\underline{A}
\label{kroneckerMatrix2}
\end{align}
The equation \eqref{kroneckerMatrix2} implies that the entries $\mathcal{A}^{l}\overline{\mathbf{p}}$ can be computed recursively, involving Lyapunov-like operators\footnote{This is not a Lyapunov equation because $\underline{P}$ is known.} of the form \eqref{kroneckerMatrix2}. That is, the a priori pattern, defined by \eqref{KroneckerRepresentation22}, can be efficiently computed by performing operations on $(Nn)\times (Nn)$ sparse (multi) banded matrices. Our numerical experience shows that the complexity of computing a priori pattern for small values of $z_{1}$, is negligible compared to the complexity of the methods for approximating $\underline{X}$, presented in Section \ref{methodsApproximation}.
\par
Finally, it should be observed that the initial guess of the a priori pattern of $\underline{X}$, given by \eqref{aprioriSeriesExpansion}, can be used for other classes of sparse matrices $\underline{A}$ and $\underline{P}$, such as non-symmetric banded matrices for example. 
\section{Methods for computing sparse approximation to $\underline{X}$}
\label{methodsApproximation}
Using the insights from previous sections, in this section we develop two methods for computing (sparse) banded approximations to $\underline{X}$.
\subsection{First method}
For presentation clarity and brevity, in this section we will restrict our attention to banded a priori patterns of $\underline{X}$. All the methods can be extended to multi-banded or even more general a priori patterns, see also Section \ref{generalSparsityPattern} and numerical experiments in Section \ref{thirdExpampleSection}. 
\par
Due to the fact that the a priori pattern is a banded matrix, the entries of $\overline{\mathbf{x}}$ corresponding to the entries of $\underline{X}$ outside the \textit{a priori bandwidth}, should be eliminated. Let $\underline{X}_{y}$ be a $y$-banded, $Nn\times Nn$, binary matrix whose non-zero entries denote the a priori pattern of $\underline{X}$. Let us assume that $y\ll nN$, that is, the matrix $\underline{X}_{y}$ is sparse and banded. Define the vector $\mathbf{x}_{y}=\text{vec}\left(\underline{X}_{y} \right)$. Let the vector $\tilde{\mathbf{x}}_{y}\in \mathbb{R}^{N_{1}}$ be defined by taking non-zero elements of $\mathbf{x}_{y}$ and stacking them on top of each other. Given that the matrix $\underline{X}_{y}$ is sparse and banded, we have that $N_{1}\ll (Nn)^{2}$ ($N_{1}$ is in the order of $Nn$). Let the matrix $\tilde{\mathcal{A}}_{1}\in \mathbb{R}^{(Nn)^{2}\times N_{1}}$ be defined by eliminating the columns of $\mathcal{A}$ corresponding to the zero entries of $\mathbf{x}_{y}$. Depending on the bandwidth $y$ and the pattern and the bandwidth of $\underline{A}$, this elimination of columns of $\mathcal{A}$ might induce zero rows of $\tilde{\mathcal{A}}_{1}$, that together with the corresponding elements of $\overline{\mathbf{p}}$ can be further eliminated. After these row eliminations of $\tilde{\mathcal{A}}_{1}$ and $\overline{\mathbf{p}}$, we obtain the matrix $\tilde{\mathcal{A}}_{2}\in \mathbb{R}^{N_{2}\times N_{1}}$ and the vector $\overline{\mathbf{p}}_{2}\in \mathbb{R}^{N_{2}}$. We assume that $N_{2}\ge N_{1}$ and that $N_{2}$ is in the order of $Nn$. The approximate solution to the Lyapunov equation can be found by solving the following least-squares problem:
\begin{align}
\min_{\tilde{\mathbf{x}}}\left\| \overline{\mathbf{p}}_{2}- \tilde{\mathcal{A}}_{2}\tilde{\mathbf{x}}\right\|_{2}^{2}
\label{KroneckerRepresentationReducedLeastSquares}
\end{align}
where $\tilde{\mathbf{x}} \in \mathbb{R}^{N_{1}}$. The solution of \eqref{KroneckerRepresentationReducedLeastSquares} can be found by solving the normal system of equations:
\begin{align}
\hat{\mathcal{A}}\tilde{\mathbf{x}}=\tilde{\mathcal{A}}^{T}_{2}\overline{\mathbf{p}}_{2}
\label{normalSystem}
\end{align}
where $\hat{\mathcal{A}}=\tilde{\mathcal{A}}^{T}_{2}\tilde{\mathcal{A}}_{2}$ and $\hat{\mathcal{A}}\in \mathbb{R}^{N_{1}\times N_{1}}$. The problem \eqref{KroneckerRepresentationReducedLeastSquares}-\eqref{normalSystem} can be efficiently solved using the Conjugate Gradient Least-Squares (CGLS) method \cite{saad2003iterative,bjorck1996numerical}, without the need to explicitly form the normal equations. One iteration of the CGLS method takes about $2\text{nz}(\tilde{\mathcal{A}}_{2})+3N_{1}+2N_{2}$ flops, where  $\text{nz}(\tilde{\mathcal{A}}_{2})$ stands for the number of non-zero elements of $\tilde{\mathcal{A}}_{2}$ \cite{bjorck1996numerical}. 
\par
It is well-known that the convergence of the CG methods is fast for well-conditioned problems\cite{saad2003iterative}. More precisely, the convergence of the CGLS method is determined by the factor\footnote{Factor closer to one implies slower convergence, and closer to zero means faster convergence.} $(\kappa_{1}-1)/(\kappa_{1}+1)$, where $\kappa_{1}$ is the condition number of $\hat{\mathcal{A}}$, for more details see Chapter 7 of \cite{bjorck1996numerical}. Due to the fact that the matrix $\tilde{\mathcal{A}}_{2}$ is formed from the matrix $\mathcal{A}$ and $\hat{\mathcal{A}}=\tilde{\mathcal{A}}_{2}^{T}\tilde{\mathcal{A}}_{2}$, it follows that the condition number of $\hat{\mathcal{A}}$ is related to the condition number of $\mathcal{A}$ or to the condition number of $\underline{A}$. However, a theoretical study of this relationship is nontrivial and it is left for further research. Our numerical simulations indicate that if the matrix $\mathcal{A}$ is well-conditioned, the matrix $\tilde{\mathcal{A}}_{2}$ inherits this numerical property. Furthermore, the matrix $\hat{\mathcal{A}}$ is also relatively well-conditioned (although the condition number of $\hat{\mathcal{A}}$ is usually larger than the condition number of $\mathcal{A}$). On the other hand, we showed that for a well-conditioned $\underline{A}$, the entries of $\underline{X}$ decay quickly away from a banded pattern or from the main diagonal, and consequently, the bandwidth $y$ of the a priori pattern $\underline{X}_{y}$ can be chosen to be much smaller than $Nn$, without seriously compromising the accuracy. That is, for a well-conditioned $\underline{A}$, we have that $N_{1}$ and $N_{2}$ are naturally in the order of $Nn$. All these observations, together with the experience gathered by performing numerical simulations, indicate that for well-conditioned problems, the solution to \eqref{normalSystem} can be determined efficiently, with $O(N)$ computational and memory complexities, see also Remark \ref{sparsePreconditioners}. Theoretical analysis that can support our observations is left for a future research.
\begin{rem}
For a relatively ill-conditioned $\underline{A}$, the convergence rate of the CGLS method can be improved by employing the preconditioning techniques. To preserve the sparsity of the problem, techniques that employ sparse approximate inverse preconditioners can be used \cite{chow1998approximate,chow2000priori,benzi1998sparse}.
\label{sparsePreconditioners}
\end{rem}

\subsection{Second method}
\label{secondMethodSection}
Here we will develop a second method for approximating $\underline{X}$. This method is based on the approximation of $\text{exp}\left(t\underline{A}\right)$ by a banded matrix. Using this approximation we approximate the integral representation \eqref{integralRepresentationLyapunov}. Such a solution is then used as an initial guess for the gradient projection method. 
\par
In \cite{benzi2007,benzi2015decay} it has been shown that the matrix exponential of a symmetric banded matrix is an off-diagonally decaying matrix. A relatively non-conservative upper bound on the off-diagonal decay rate of the matrix exponential is derived in Theorem 4.2 in \cite{benzi2015decay}. This decay bound depends on the extreme eigenvalues $a$ and $b$ of $\underline{A}$. However, the influence of the condition number on the decay rate of $\text{exp}\left(t\underline{A}\right)$, is still an open problem\footnote{It should be expected that the decay rate is faster for well-conditioned matrices, because it very well known that the decay rate of matrix functions, such as $\underline{A}^{-1}$, is faster if the condition number of $\underline{A}$ is smaller.}. Due to the fact that $\exp\left(t\underline{A} \right)$ is an off-diagonally decaying matrix, it can be approximated by a banded matrix using the Chebyshev series \cite{benzi2007,bergamaschi2000efficient}, see also Remark \ref{nonSymChebyshev}. First, we transform the matrix function $\exp\left(t\underline{A}\right)$ into a complex function $\exp\left(tz\right)$, where $z$ is a complex number, belonging to a domain that contains the spectrum of $\underline{A}$, that is, $z\in [a,b]$. It is obvious that the eigenvalues of $t\underline{A}$ are in the interval $[ta,tb]$. The next step is to transform this interval into the interval $[-1,1]$. This can be achieved by defining a new variable $w$ as follows:
\begin{align}
w=\frac{2tz-t(a+b)}{t(b-a)}
\label{shiftAndScale}
\end{align}
It can be easily seen that when $tz\in [ta,tb]$ then $w \in [-1,1]$. From \eqref{shiftAndScale} we have:
\begin{align}
tz=\frac{t}{2}\left( \left(b-a\right)w+a+b\right)
\label{shiftAndScaleExpressed}
\end{align}
Using \eqref{shiftAndScaleExpressed}, we have:
\begin{align}
\exp\left(tz\right)=\exp\left(\frac{t}{2}\left( \left(b-a\right)w+a+b \right)\right)
\label{eMatrixExpdecay}
\end{align}
Let $\underline{A}_{1}$ be a matrix corresponding to the complex variable $w$. By substituting $w$ with $\underline{A}_{1}$ and $z$ by $\underline{A}$ in \eqref{shiftAndScale}, we obtain:
\begin{align}
\underline{A}_{1}=\frac{2}{b-a}\underline{A}-\frac{a+b}{b-a}I
\label{matrixFunctions1}
\end{align}
The eigenvalues of $\underline{A}_{1}$ belong to the interval $[-1,1]$. Similarly, from \eqref{shiftAndScaleExpressed} we obtain:
\begin{align}
t\underline{A}=\frac{t}{2}\left( (b-a)\underline{A}_{1}+(a+b) I \right) \label{matrixFunctions2} 
\end{align}
By substituting $t\underline{A}$ in $\exp\left(t\underline{A}\right)$ with \eqref{matrixFunctions2}, we define the function $f\left(\underline{A}_{1}\right)$:
\begin{align}
f\left(\underline{A}_{1}\right)=\exp\left(\frac{t}{2}\left( (b-a)\underline{A}_{1}+(a+b) I \right)\right)
\label{functionDefinition}
\end{align}
It is obvious that $f\left(\underline{A}_{1}\right)=\exp\left(t\underline{A}\right)$. Consequently, the behavior of $f\left(\underline{A}_{1}\right)$ is identical to the behavior of  $\exp\left(t\underline{A}\right)$.
\par
 Consider the matrix $\underline{A}_{1}$ defined in \eqref{matrixFunctions1}. Taking into account that the spectrum of $\underline{A}_{1}$ belongs to the interval $[-1,1]$, the truncated Chebyshev series expansion of the matrix exponential \eqref{functionDefinition} is defined by \cite{benzi2007,mason2003chebyshev}:
\begin{align}
&f(\underline{A}_{1})\approx \tilde{f}(\underline{A}_{1}),\notag \\
&\tilde{f}(\underline{A}_{1})=\frac{c_{1}}{2}I+\sum_{k=2}^{M}c_{k}T_{k}\left( \underline{A}_{1} \right)=\sum_{k=1}^{M}{}^{\prime} c_{k}T_{k}\left( \underline{A}_{1} \right)
\label{truncatedChebyshevSeries}
\end{align}
where $c_{k} \in \mathbb{R}$ are the Chebyshev coefficients, $T_{k}\left(\underline{A}_{1}\right)\in \mathbb{R}^{Nn\times Nn}$ are the Chebyshev (matrix) polynomials of the first kind, and the symbol $\sum_{k=1}^{M}{}^{\prime}$ means that the first term in the sum is halved \cite{Mathar2006}. The Chebyshev matrix polynomials are defined by \cite{benzi2007,mason2003chebyshev}:
\begin{align}
T_{1}=I,\; T_{2}=\underline{A}_{1}, \; T_{k+1}=2\underline{A}_{1}T_{k}-T_{k-1},\; k=2,3,\ldots
\label{ChebyshevPolynomialsDefinition}
\end{align}
The Chebyshev coefficients can be approximated by \cite{benzi2007,mason2003chebyshev} (see Remark \ref{remarkChebyCoefficient}):
\begin{align}
c_{k}\approx \frac{2}{R} \sum_{j=1}^{R}f\left(\cos\left( \theta_{j}\right) \right)\cos\left( \left(k-1 \right)\theta_{j}\right)
\label{theChebyConstants}
\end{align}
where $\theta_{j}=\pi \left(j-\frac{1}{2} \right)/R$ and $R$ is a sufficiently large positive integer. The function $f\left(\cos\left( \theta_{j}\right)\right)$ in \eqref{theChebyConstants} is defined as follows. First, by formally substituting in \eqref{functionDefinition} the matrix $\underline{A}_{1}$ with $w$, and $I$ with $1$, we define:
\begin{align}
f\left(w\right)=\exp\left(\frac{t}{2}\left( (b-a)w+a+b \right)\right)
\label{temp111}
\end{align}
Then, substituting in \eqref{temp111} the argument $w$ with $\cos\left( \theta_{j}\right)$ we obtain:
\begin{align}
f\left(\cos\left( \theta_{j}\right) \right)=\exp\left(t\hat{q}_{j}\right), \; \hat{q}_{j}=\frac{ (b-a)\cos\left( \theta_{j}\right)+a+b}{2}
\notag
\end{align}
\begin{rem}
An alternative method for computing the Chebyshev coefficients of the exponential function is based on the Bessel functions, see the equation (2.1) in \cite{bergamaschi2000efficient}. $\hfill \square$
\label{remarkChebyCoefficient}
\end{rem}
\begin{rem}
In the case of non-symmetric matrices $\underline{A}$, the Chebyshev approximation method cannot be directly used to approximate the matrix exponential. In the case of non-symmetric problems, the Faber polynomials need to be used, as demonstrated in \cite{bergamaschi2003efficient}. 
\label{nonSymChebyshev}
\end{rem}
The Chebyshev approximation error is defined by $\epsilon_{M}=\left\| f(\underline{A}_{1})- \tilde{f}(\underline{A}_{1})\right\|_{2}$. It can be easily shown that (see for example, Section 4.1 in \cite{benzi2007}):
\begin{align}
\epsilon_{M}\le \sum_{k=M+1}^{\infty} |c_{k}|
\label{ChebApproxErrorUpper}
\end{align}
That is, the approximation error does not depend (explicitly) on the dimensions of  $\underline{A}_{1}$. On the other hand, in \cite{bergamaschi2000efficient} it has been shown that:
\begin{align}
\epsilon_{M}\le \sum_{k=M+1}^{\infty} |c_{k}|\le 2\exp\left(tb \right)K_{2}(M,a,t)
\label{ChebApproxErrorUpper22}
\end{align}
where $K_{2}(M,a,t)$ is a constant depending on $M$, $b$ and $t$, for more details see the equations (2.2) and (2.3) in \cite{bergamaschi2000efficient}. By approximately computing the right-hand side of \eqref{ChebApproxErrorUpper}, or by computing \eqref{ChebApproxErrorUpper22}, we can find the maximum order of the Chebyshev polynomials $M$ for which the approximation error is below a predefined small number. However, if $M$ is large, then the Chebyshev approximation $\tilde{f}(\underline{A}_{1})$ is a dense matrix. Namely, from \eqref{ChebyshevPolynomialsDefinition} we see that each Chebyshev polynomial can be expressed as the sum of powers of $\underline{A}_{1}$. Due to the fact that $\underline{A}_{1}$ is $m$-banded, $\underline{A}_{1}^{k}$ is $km$-banded, $k=1,\ldots, M$. That is, for large $k$, the matrix $\underline{A}_{1}^{k}$ becomes fully populated. Given that $N$ is large, it might be impossible to compute and to store $\underline{A}_{1}^{k}$. This implies that in order to ensure that the matrix $\tilde{f}(\underline{A}_{1})$ is (sparse) banded and to guarantee that it can be computed and stored with $O(N)$ complexity, $M$ needs to be kept small. If the off-diagonal decay rate of $\exp\left(t\underline{A} \right)$ is fast, then our numerical results show that even for a small $M$, the approximation accuracy is relatively good. However, if the off-diagonal decay rate is not fast, selecting a relatively small $M$ can seriously compromise the approximation accuracy. 
\\
One of the ways to resolve this problem is to employ the \textit{numerical dropping}  technique \cite{benzi2007,haber2014Gramian}. Namely we introduce the following modification in \eqref{ChebyshevPolynomialsDefinition}:
\begin{align}
 T_{k+1}=\mathcal{D}\left(2\underline{A}_{1}T_{k}-T_{k-1}\right),\; k=2,3,\ldots
\label{sparsifiedChebyshev}
\end{align}
where the \textit{bandwidth projection operator} $\mathcal{D}\left( Z \right)$, acting on an arbitrary matrix $Z=[z_{i,j}]$, is defined by:
\begin{align}
\mathcal{D}\left( Z \right)=\left\{ \begin{array} {rl}
z_{i,j}, &   \;\;  |i-j|\le d/2 \\
0, &  \;\; |i-j|> d/2 \end{array} \right. 
\label{sparsificationOperatordefinition}
\end{align}
where $d$ is the prescribed bandwidth. By setting to zero entries of $2\underline{A}_{1}T_{k}-T_{k-1}$ that are outside the bandwidth $d$ (that should be relatively small compared to $Nn$), the bandwidth projection operator ensures that each Chebyshev polynomial $T_{k+1}$ remains a (sparse) banded matrix. In this way we can select a large $M$ and ensure that the matrix $\tilde{f}(\underline{A}_{1})$ can be computed with $O(N)$ memory and computational complexities, for more details see Section \ref{complexityAnalysis}. The parameter $d$ can be selected using the results of Theorem 4.2 in \cite{benzi2015decay}, or using simple heuristics. Moreover, using the 
results of \cite{simoncini2015lyapunov,canuto2014decay,benzi2015decay} it is possible to  construct a projection operator in \eqref{sparsifiedChebyshev}, that projects $T_{k+1}$ onto a more optimal matrix pattern. Finally, the results of Section \ref{generalSparsityPattern} can be extended to construct a more optimal projection operator.
\\
We have observed that for small values of $d$ and for relatively ill-conditioned $\underline{A}$, the modified Chebyshev recurrence \eqref{sparsifiedChebyshev} might diverge. To improve the accuracy of the Chebyshev approximation and to increase the stability of \eqref{sparsifiedChebyshev}, instead of applying $\mathcal{D}\left(\cdot\right)$ in every iteration of \eqref{sparsifiedChebyshev}, it can be applied in every second or third iteration. However, this strategy increases the computational complexity.
\par
After we presented the method for approximating $\text{exp}\left(t\underline{A}\right)$, we turn our attention to the problem of approximating the integral in \eqref{integralRepresentationLyapunov}. From Theorem 4.3 in \cite{grasedyck2003solution}, we have that $\underline{X}$ can be approximated by:
\begin{align}
\tilde{\underline{X}}=-\sum_{j=-q}^{q} \psi \omega_{j}\exp\left( \psi t_{j} \underline{A} \right)\underline{P}\exp\left( \psi t_{j} \underline{A} \right)
\label{method1LyapunovApproximate}
\end{align}
where $q$ is a positive integer, and
\begin{align}
&\psi=\frac{3}{2|b+\epsilon_{1}|}, \; \; \omega_{j}=\left(q+q\exp\left(-2jq^{-1/2} \right)\right)^{-1/2} \notag \\
&t_{j}=\log\left( \exp\left(jq^{-1/2} \right)+\sqrt{1+\exp\left(2jq^{-1/2}\right)} \right) 
\label{method1LyapunovApproximateExplanation}
\end{align}
where $\epsilon_{1}\ll |b|$ is a small number. In \cite{grasedyck2003solution} it has been shown that the approximation error exponentially decreases with $\sqrt{q}$:
\begin{align}
\left\|\underline{X}-\tilde{\underline{X}}\right\|_{2}\le K\left(\underline{A} \right)\left\|\underline{P} \right\|_{2}\exp\left(-\sqrt{q} \right)
\label{exponentialDecreaseFirstMethod}
\end{align}
where the constant $K\left(\underline{A} \right)$ depends on $a$ and $b$. That is, as $q$ approaches infinity, the approximation error approaches zero.
\\
Consider the matrix $\exp\left( \psi t_{j} \underline{A} \right)$ in \eqref{method1LyapunovApproximate}. Let $\tilde{t}_{j}=\psi t_{j}$. For each $\tilde{t}_{j}$ we can compute the Chebyshev approximation  $\tilde{f}_{j}\left( \underline{A}_{1}\right)$ of the matrix exponential $\exp\left( \tilde{t}_{j} \underline{A} \right) $. By substituting in \eqref{method1LyapunovApproximate} the matrix $\exp\left( \tilde{t}_{j} \underline{A} \right) $ with $\tilde{f}_{j}\left( \underline{A}_{1}\right)$, we define the approximate solution of the Lyapunov equation:
\begin{align}
\tilde{\underline{X}}_{1}=-\sum_{j=-q}^{q} \psi \omega_{j}\tilde{f}_{j}\left(\underline{A}_{1}\right)\underline{P}\tilde{f}_{j}\left(\underline{A}_{1}\right)
\label{method1LyapunovApproximateFinal}
\end{align}
Let us assume that each of the matrices $\tilde{f}_{j}\left( \underline{A}_{1}\right)$, $j=-q,\ldots,q$, are computed using the Chebyshev approximation with the bandwidth projection operator \eqref{sparsifiedChebyshev}. Consequently, the matrices $\tilde{f}_{j}\left( \underline{A}_{1}\right)$ are sparse, banded matrices with the bandwidth equal to $d$. Given that the bandwidth of $\underline{P}$ is equal to $l$, we have that the total bandwidth of $\tilde{f}_{j}\left( \underline{A}_{1}\right)\underline{P}\tilde{f}_{j}\left(  \underline{A}_{1}\right)$ is $2d+l$. Given that $d\ll N$ and $l\ll N$, we have that $2d+l\ll N$, that is, the matrix $\tilde{f}_{j}\left(\underline{A}_{1}\right)\underline{P}\tilde{f}_{j}\left( \underline{A}_{1}\right)$ is a sparse banded matrix. Furthermore, because the sum of matrices of equal bandwidths does not increase the bandwidth of the resulting sum, we have that $\tilde{\underline{X}}_{1}$ is a sparse banded matrix with the bandwidth equal to $2d+l$. If $q\ll N$, then it is obvious that \eqref{method1LyapunovApproximateFinal} can be computed with $O(N)$ computational and memory complexities, for more details see Section \ref{complexityAnalysis}.
\par
Next, we improve the accuracy of the approximate solution $\tilde{\underline{X}}_{1}$ by using it as an initial guess of the decision variable of the constrained matrix least-squares problem \cite{li2012weighted}:
\begin{align}
&\min_{\underline{X}}\left\|\underline{P}-\underline{A}\underline{X}-\underline{X}\underline{A}^{T}\right\|_{F}^{2} 
\label{solutionLyapunovOptimization} \\
& \text{subject to}\;\; \underline{X} \in \mathcal{X}_{d_{1}} \label{solutionLyapunovOptimizationConstraint}
\end{align}
where $\mathcal{X}_{d_{1}}$ denotes the set of all banded matrices with the bandwidth of $d_{1}$. The \textit{gradient projection method} for solving \eqref{solutionLyapunovOptimizationConstraint} has the following form \cite{bertsekas1999nonlinear}:
\begin{align}
&\underline{X}_{k+1}=\mathcal{D}_{1}\left(\underline{X}_{k}-\delta_{k} \mathcal{G}_{k}\right),\; k=0,1,2,\ldots
\label{steepestDescentMethod} 
\end{align}
where $\mathcal{D}_{1}\left( \cdot \right)$ is the bandwidth projection operator defined in \eqref{sparsificationOperatordefinition} for the bandwidth of $d_{1}$, $k$ is the iteration index, $\delta_{k}$ is the step size and $\mathcal{G}_{k}$ is the gradient defined by \cite{li2012weighted}:
\begin{align}
&\mathcal{G}_{k} =-2\underline{A}^{T}\mathcal{R}_{k}-2\mathcal{R}_{k}\underline{A}
\label{steepestDescentExplicit} \\
&\mathcal{R}_{k}=\underline{P}-\underline{A}\underline{X}_{k}-\underline{X}_{k}\underline{A}^{T} 
\label{steepestDescentExplicitExplanation}
\end{align}
To define the step size $\delta_{k}$, we first define the following two quantities:
\begin{align}
& F_{1}\left(\underline{X}_{k}\right)=\left\|\underline{P}-\underline{A}\underline{X}_{k}-\underline{X}_{k}\underline{A}^{T}\right\|_{F}^{2}, \notag \\ &\underline{X}_{k}\left( \delta \right)=\mathcal{D}_{1}\left(\underline{X}_{k}-\delta \mathcal{G}_{k}\right)
\end{align}
Keeping these definitions in mind, the step size is determined by the Armijo rule along the projection arc \cite{bertsekas1999nonlinear}:
\begin{align}
\delta_{k}=\zeta^{h_{k}} \overline{\delta}
\label{stepSize}
\end{align}
where $h_{k}$ is equal to the first nonnegative integer $h$ for which:
\begin{align}
F_{1}\left(\underline{X}_{k}\right)-F_{1}\left(\underline{X}_{k}\left( \zeta^{h} \overline{\delta}\right)\right)\ge \sigma \mathcal{G}_{k}^{T}\left(\underline{X}_{k}-\underline{X}_{k}\left(\zeta^{h} \overline{\delta} \right) \right)
\label{ArmijoRule}
\end{align}
and $\sigma \in (0,1)$, $\zeta \in (0,1)$ and $\overline{\delta}>0$.
\textit{To summarize, the proposed method consists of the following two steps
\begin{enumerate}
\item Compute the matrix $\tilde{\underline{X}}_{1}$ defined in \eqref{method1LyapunovApproximateFinal}.
\item Set $\underline{X}_{0}=\tilde{\underline{X}}_{1}$ and propagate the iteration  \eqref{steepestDescentMethod} until convergence or the maximal number of iterations has been reached. 
\end{enumerate}}
The convergence rate of the projected gradient method is well-studied in the literature, see for example \cite{bertsekas1999nonlinear}. Briefly speaking, the convergence rate is  similar to the convergence rate of the unconstrained steepest descent method, and it depends on the condition number of $\underline{A}$. For well-conditioned problems the convergence is relatively fast. However, for ill-conditioned problems it can be slow. One of the ways to improve the convergence for ill-conditioned problems is to use scaled gradient projection methods, for more details see Chapter 2 of \cite{bertsekas1999nonlinear}. The optimization problem \eqref{solutionLyapunovOptimization}-\eqref{solutionLyapunovOptimizationConstraint}
can also be solved using other methods, for example using the matrix version of MINRES method with numerical droppings \cite{saad2003iterative,benner2013efficient,chow1998approximate}.

\subsection{Complexity analysis}
\label{complexityAnalysis}
\par
Let us analyze the complexity of the steps needed to compute the initial guess. The first computationally demanding step is to compute the extreme eigenvalues $a$ and $b$ of $\underline{A}$. This step can be performed efficiently using the ARPACK software or MATLAB functions $\text{eigs}\left(\cdot\right)$ and  $\text{svds}\left(\cdot\right)$, that are based on the Implicitly Restarted Arnoldi methods \cite{lehoucq1998arpack}. Generally speaking, one iteration of the methods implemented in ARPACK  or in MATLAB, can be computed with $O(N)$ complexity, where the convergence rate depends on the spectral properties of $\underline{A}$. The storage requirement for computing the extreme eigenvalues is $O(N)$ \cite{lehoucq1998arpack}. Due to the fact that the initial guess is later on improved using the gradient projection method, we can tolerate some inaccuracies in computing extreme eigenvalues. That is, the tolerances for computing the extreme eigenvalues can be relaxed with the goal of achieving overall $O(N)$ computational complexity. 
\par
Let us now turn our attention to the Chebyshev approximation of the matrix exponential \eqref{truncatedChebyshevSeries} and the integral approximation \eqref{method1LyapunovApproximateFinal}. In Appendix \ref{complexitySecond}, we have shown that the computational complexity of these steps scales approximately linearly with $N$, $M$ and $q$ (assuming all other parameters are fixed), while on the other hand, it scales quadratically with $d$, where it is assumed that the bandwidth projection operator \eqref{sparsifiedChebyshev} is applied. For a relatively precise operations count see Appendix \ref{complexitySecond}. Memory complexity scales with $O(N)$, assuming that $d\ll N$. Similarly it can be shown that the computational complexity of one iteration of the gradient projection method \eqref{steepestDescentMethod} scales linearly with $N$ and $d_{1}$, while the memory complexity scales with $O(N)$, assuming that $d_{1}\ll N$.
\section{Numerical Experiments}
\label{sectionNumericalExperiments}
We demonstrate the effectiveness of the developed methods on three examples. Numerical simulations are performed in MATLAB on a laptop computer with $4$ GB of RAM, with processor specifications: Intel(R) Core(TM) i5-2410M CPU @ $2.30$ GHz  $2.30$ GHz.
\subsection{First example: 2D heat equation}
\begin{figure}[H]
  \centering
 \includegraphics[scale=0.21,trim=0mm 0mm 0mm 0mm ,clip=true]{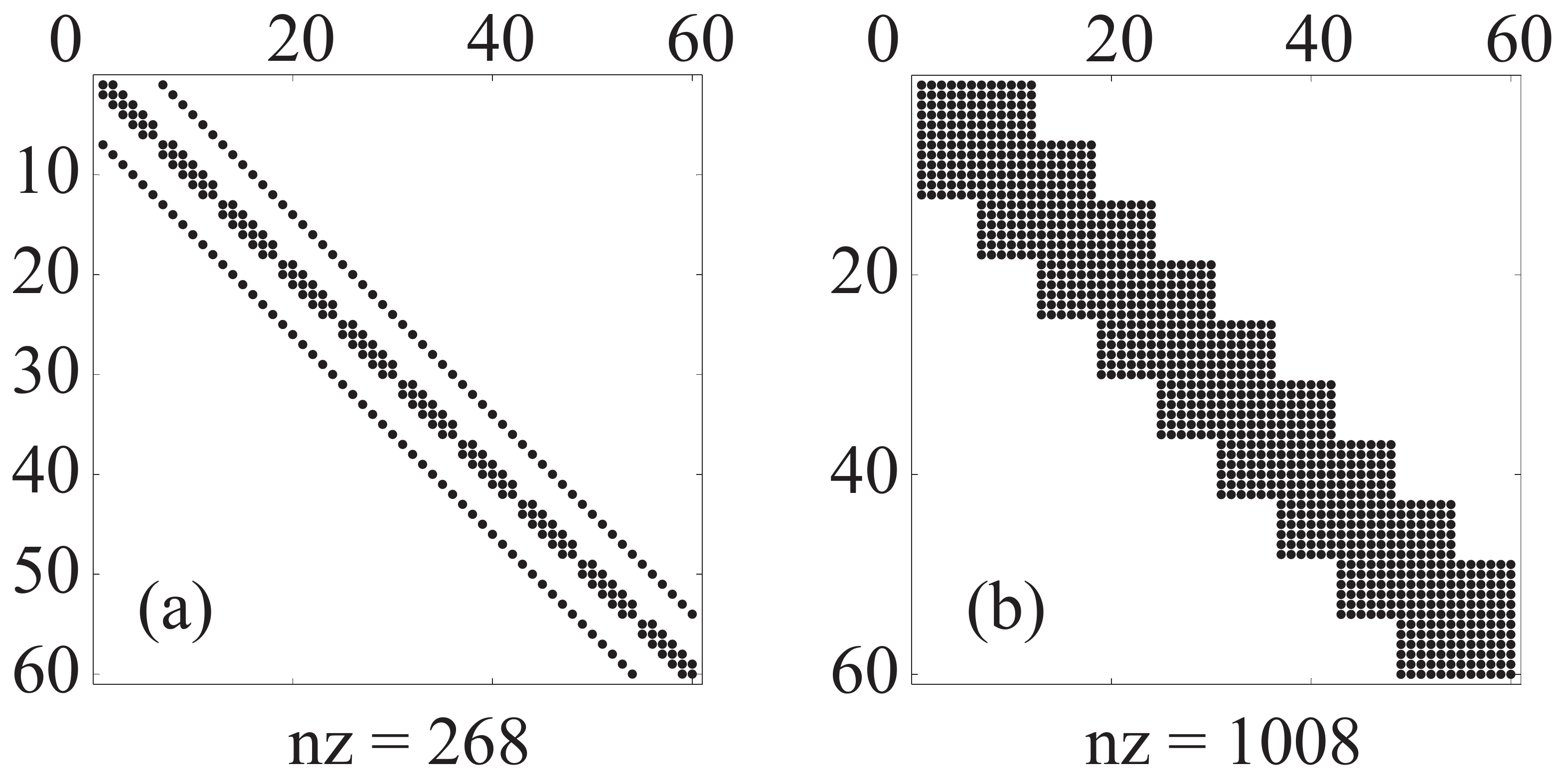}
 \caption{Sparsity patterns: (a) $\underline{A}$. (b) $\underline{P}$, "nz" denotes the number of non-zero elements. Results are generated for $N=10$.}
\label{fig:SparistyPatterns}
\end{figure}

\begin{figure}[H]
  \centering
 \includegraphics[scale=0.23,trim=0mm 0mm 0mm 0mm ,clip=true]{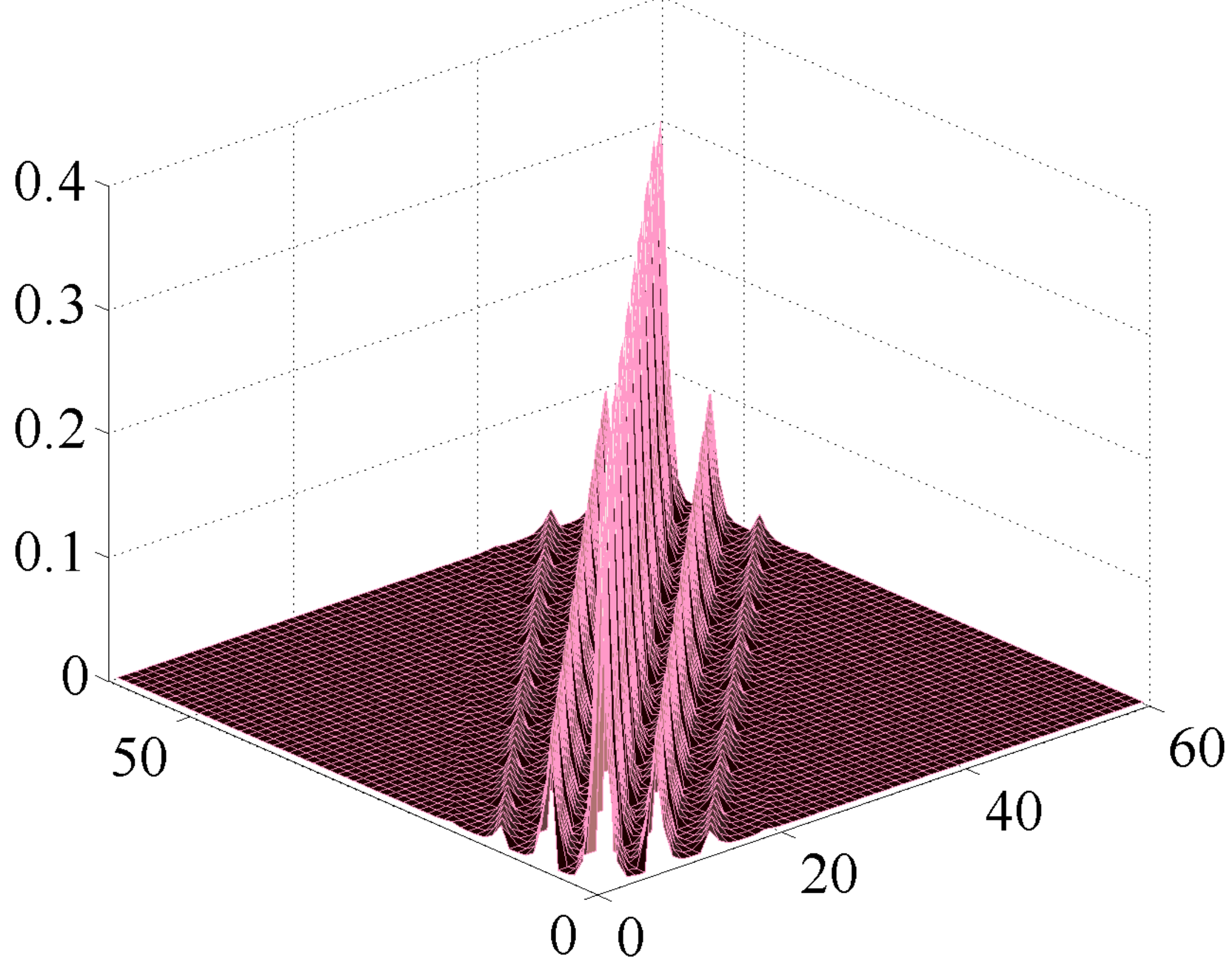}
   \caption{Surface plot of $\exp\left(t\underline{A} \right)$ for $t=1$.}
\label{fig:surfPlotExp}
\end{figure}

\begin{figure}[H]
  \centering
 \includegraphics[scale=0.2,trim=0mm 0mm 0mm 0mm ,clip=true]{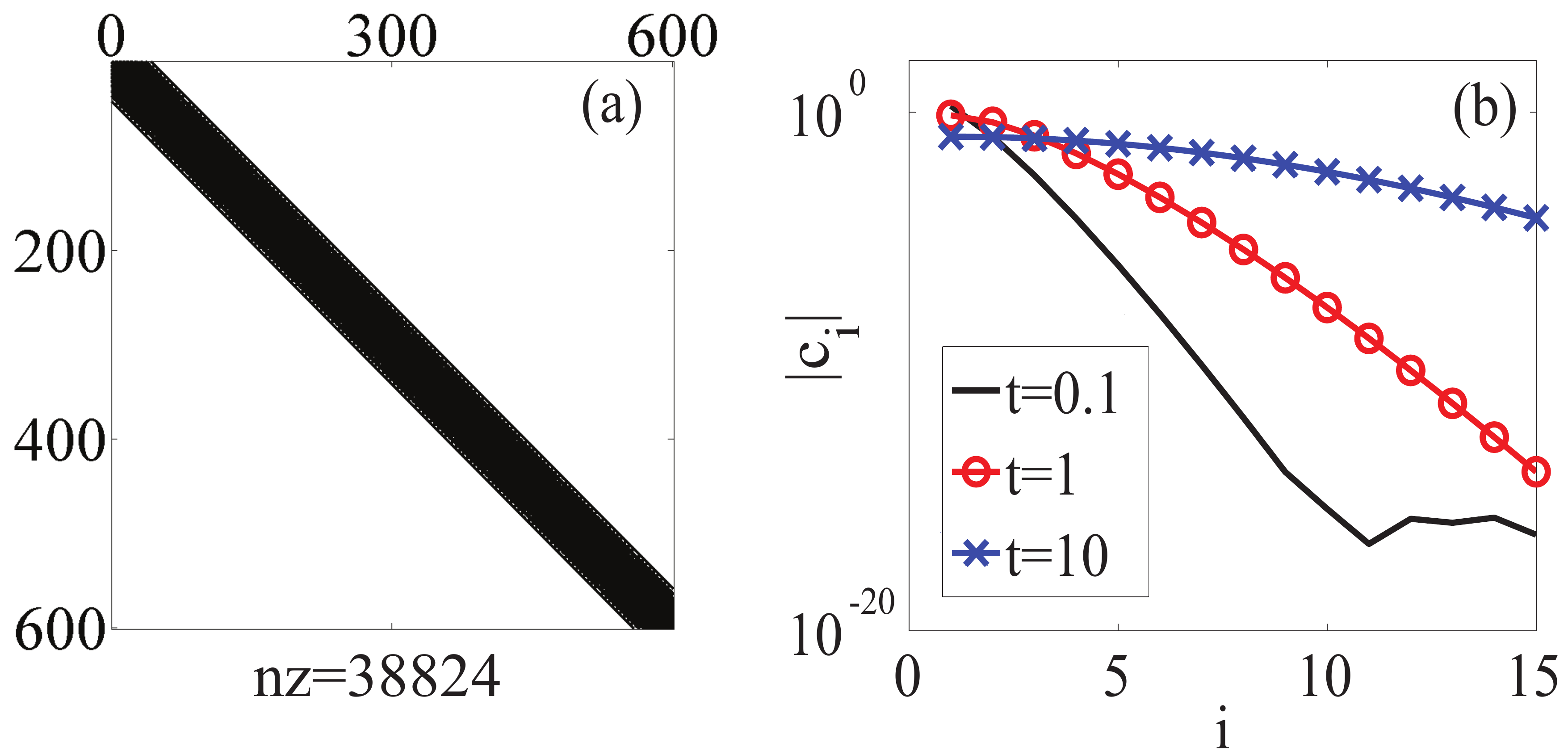}
   \caption{(a) Sparsity pattern of the approximation of $\exp\left(\underline{A} \right)$. (b) The dependence of the Chebyshev coefficients on $t$.}
\label{fig:sparsityExpApp}
\end{figure}

We consider a model describing the temperature change of a thermally actuated deformable mirror used in extreme ultraviolet lithography \cite{haberThesis,haber2013OL,haber2013predictive,Haber:13}. Heat equation constants, discretization grid and discretization steps are defined in Chapter 2 of \cite{haberThesis}.
The structure of this model is equivalent to a model obtained by the finite-difference discretization of the Laplace operator on a rectangular domain. The matrices $\underline{A}$ and $\underline{P}$ are block tri-diagonal, with the main block diagonals defined by:
\vspace{-6mm}
\begin{small}
\begin{align}
&A_{i,i}=\begin{bmatrix}a & e & 0 & 0 & 0 & 0\\ e & a & e & 0 & 0 & 0 \\ 0 & e & a & e & 0 & 0\\ 0 & 0 & e & a & e & 0 \\  0 & 0 & 0 & e & a & e \\ 0 & 0 & 0 & 0 & e & a  \end{bmatrix},\; 
P_{i,i}=-\begin{bmatrix}1 & 0.2 & 0.2 & 0.2 & 0.2 & 0.2\\ 0.2 & 1 & 0.2 & 0.2 & 0.2 & 0.2 \\ 0.2 & 0.2 & 1 & 0.2 & 0.2 & 0.2\\ 0.2 & 0.2 & 0.2 & 1 & 0.2 & 0.2 \\  0.2 & 0.2 & 0.2 & 0.2 & 1 & 0.2 \\ 0.2 & 0.2 & 0.2 & 0.2 & 0.2 & 1  \end{bmatrix}
\label{discretizedHeatEquation}
\end{align}
\end{small}
The off-diagonal blocks of $\underline{A}$ are given by $A_{i,i-1}=A_{i,i+1}=eI_{6,6}$, where $a=-1.36$, $e=0.34$ and $I_{6,6}$ is a $6\times 6$ identity matrix. For the matrix $\underline{P}$ we have: $P_{i,i-1}=P_{i,i+1}=[-0.1]$ ($6\times 6$ matrices whose entries are equal to $-0.1$). 
The total number of local subsystems $N$ will be varied.
The sparsity patterns of $\underline{A}$ and $\underline{P}$, for $N=10$, are shown in Fig.\ref{fig:SparistyPatterns}.
We start with the Chebyshev approximation of the matrix exponential. For $t=1$, we compute $\exp(t\underline{A})$ using the built-in MATLAB function $\text{expm}\left(\cdot \right)$. The surface plot ("city plot") of this matrix is shown in Fig. \ref{fig:surfPlotExp}.
 Next, for $t=1$ and $N=100$, we approximate $\exp(t\underline{A})$ using the Chebyshev method. For $M=7$ in \eqref{truncatedChebyshevSeries}, the approximation error is $\epsilon_{M}=4.4\times 10^{-7}$, and the sparsity pattern of the approximate matrix exponential is shown in Fig. \ref{fig:sparsityExpApp}(a). These results confirm that $\exp\left(\underline{A} \right)$ can be approximated by a sparse banded matrix with high accuracy \cite{benzi1999bounds,benzi2007}.  
In Fig. \ref{fig:sparsityExpApp}(b), we illustrate the dependence of the Chebyshev coefficients on time $t$. It can be observed that as $t$ increases, the Chebyshev coefficients decay more slowly. This numerically illustrates very-well known fact that for the fixed approximation order $M$, the accuracy of approximating $\exp\left(t\underline{A} \right)$ using the Chebyshev approximation is better if $t$ is smaller \cite{benzi1999bounds,benzi2007}. This also implies that as $t_{j}$ is increased in \eqref{method1LyapunovApproximate}, to keep the accuracy of approximating 
$\exp\left(\psi t_{j}\underline{A}\right)$ constant, we need to increase the order of the Chebyshev approximation $M$.
\par
\begin{figure}[H]
  \centering
 \includegraphics[scale=0.30,trim=0mm 0mm 0mm 0mm ,clip=true]{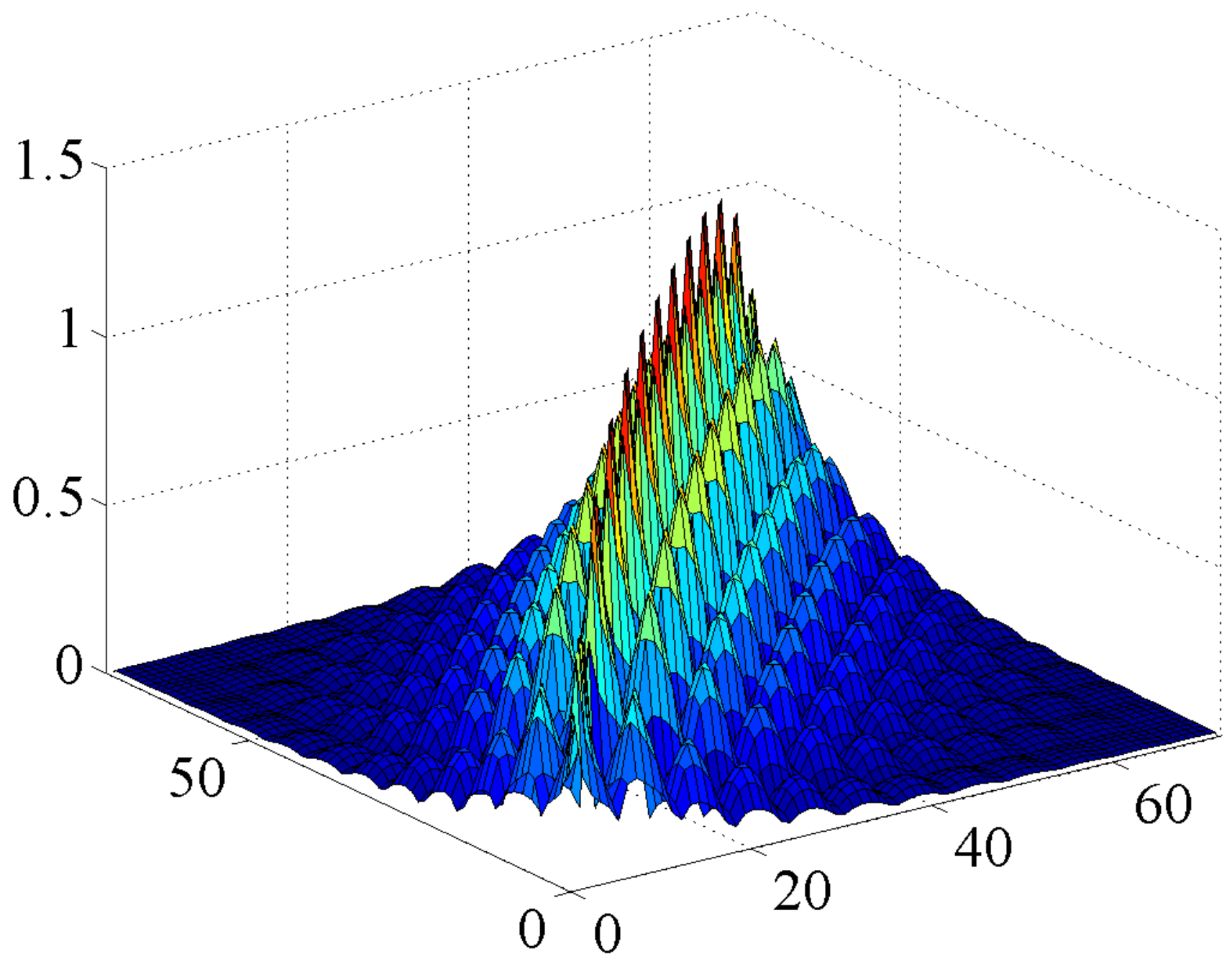}
 \caption{The surface plot of the "true" solution $\underline{X}_{T}$ computed using the function $\text{lyap}\left(\cdot \right)$. $N=12$.}
\label{fig:surfPlotXtrue}
\end{figure}
\begin{figure}[H]
  \centering
 \includegraphics[scale=0.22,trim=0mm 0mm 0mm 0mm ,clip=true]{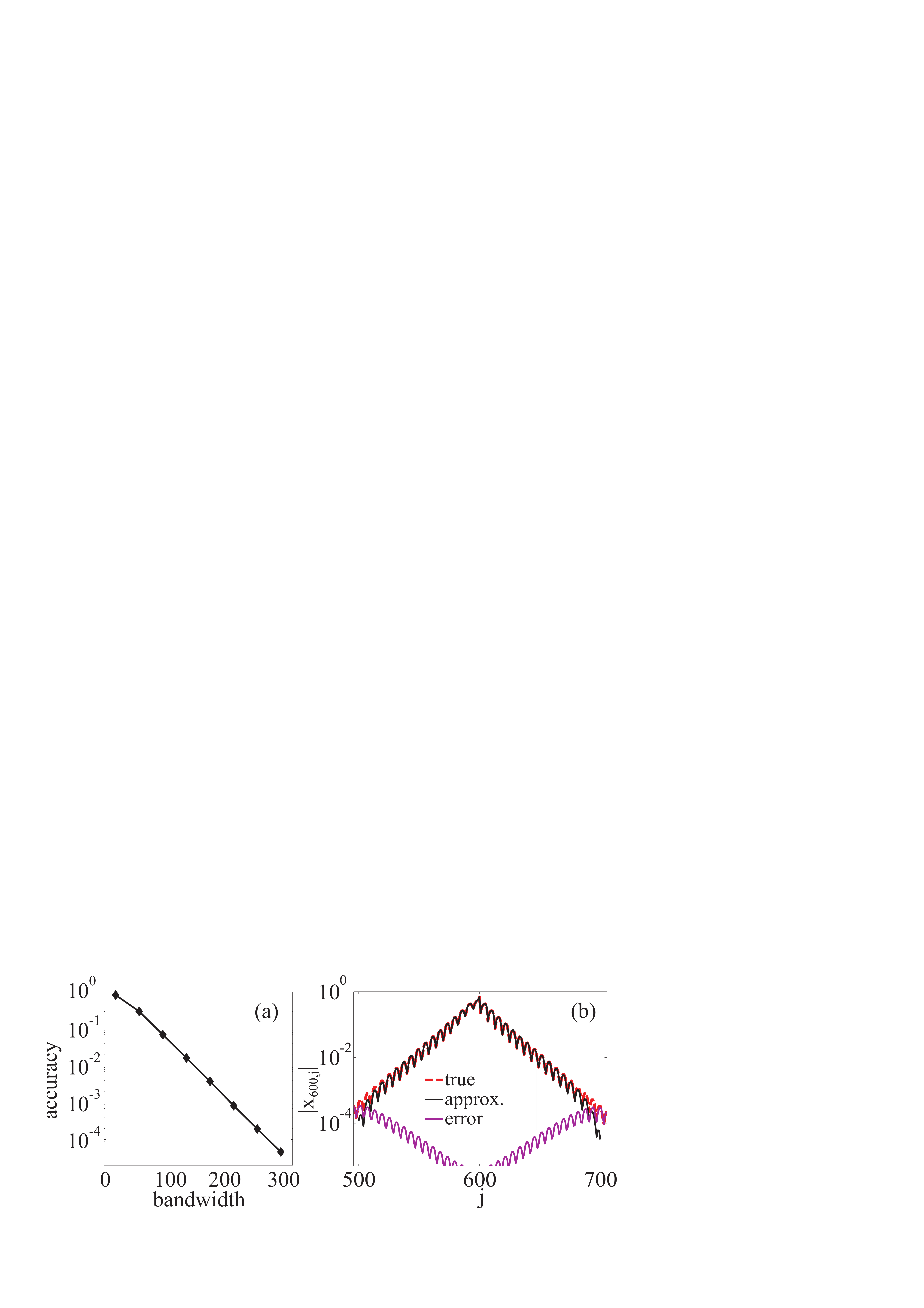}
 \caption{(a) Dependence of the accuracy \eqref{approximateSolutionAccuracy} on the bandwidth of the approximate solution. (b) Absolute value of the row $600$ of the "true" and approximate solutions, and the error. Both (a) and (b) are generated for $N=200$.}
\label{fig:Error}
\end{figure}
\begin{figure}[H]
  \centering
 \includegraphics[scale=0.23,trim=0mm 0mm 0mm 0mm ,clip=true]{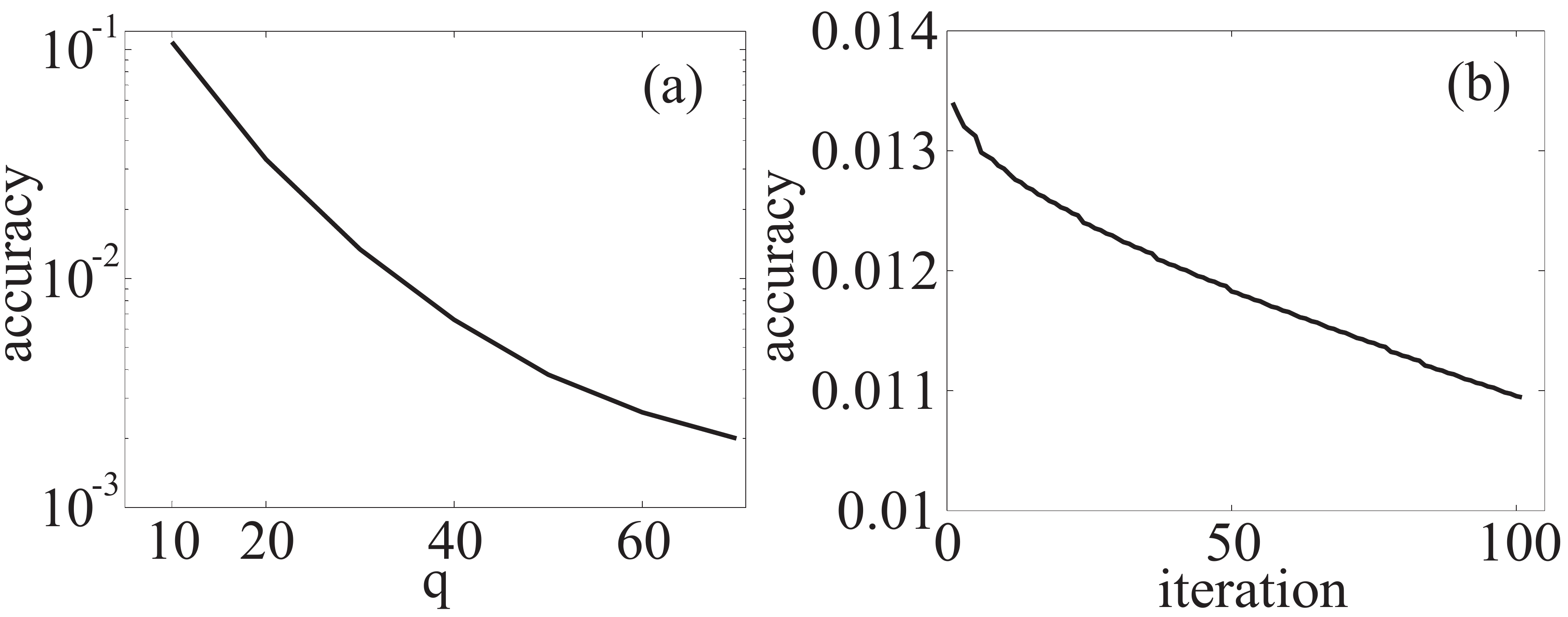}
 \caption{(a) Dependence of the accuracy \eqref{approximateSolutionAccuracy} on the parameter $q$ in the equation \eqref{method1LyapunovApproximateFinal}. (b) Accuracy dependence on the number of iterations of the gradient projection method \eqref{steepestDescentMethod}. Results are generated for $N=250$ and bandwidth of $140$.}
\label{fig:ErrorConvergence}
\end{figure}
\begin{figure}[H]
  \centering
 \includegraphics[scale=0.22,trim=0mm 0mm 0mm 0mm ,clip=true]{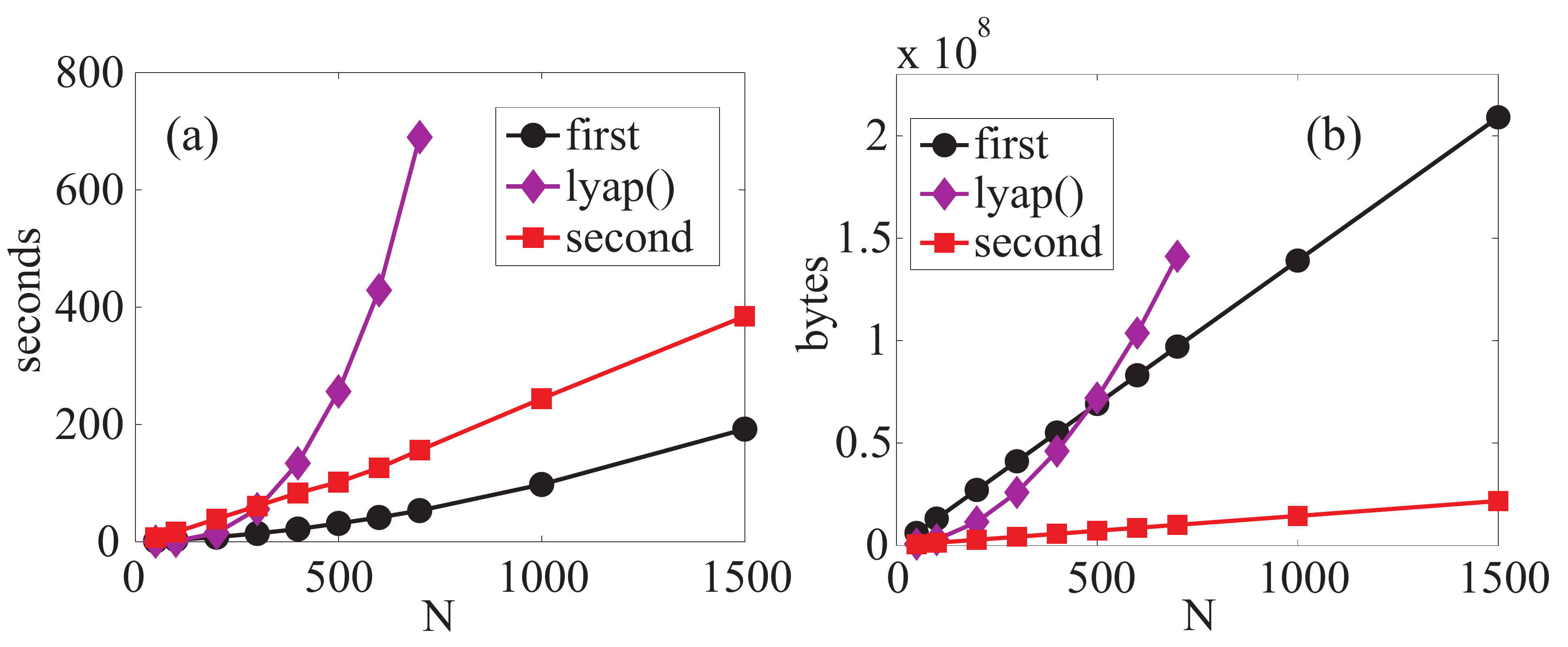}
 \caption{Complexity of the first method, the second method and the MATLAB function $\text{lyap}\left(\cdot\right)$. (a) Computational complexity. (b) Memory complexity.}
\label{fig:Complexity}
\end{figure}

We now turn our attention to the problem of approximating $\underline{X}$. We will compare the approximate solution computed using the proposed methods with the "true" solution $\underline{X}_{T}$ computed using the built-in MATLAB function $\text{lyap}\left(\cdot \right)$. This solution is a dense matrix and its surface plot is shown in Fig. \ref{fig:surfPlotXtrue}. The oscillatory behavior of entries of $\underline{X}_{T}$ is due to the Kronecker sum structure of $\underline{A}$, see Section \ref{generalSparsityPattern}.
Our first goal is to numerically investigate the accuracy, and then the complexity of the proposed methods. We first investigate how the bandwidth of the approximate solutions affects the accuracy. We only show the results for the first method because the results for the second method are similar. The stopping criteria for the CGLS method is based on the following scalar \cite{bjorck1996numerical}:
\begin{align}
\eta=\frac{\left\|  \tilde{\mathcal{A}}^{T}_{2}\left(\overline{\mathbf{p}}_{2}- \tilde{\mathcal{A}}_{2}\tilde{\mathbf{x}}^{k} \right)\right\|_{2}}{\left\| \tilde{\mathcal{A}}^{T}_{2}\left(\overline{\mathbf{p}}_{2}- \tilde{\mathcal{A}}_{2}\tilde{\mathbf{x}}^{0} \right)\right\|_{2}}
\label{stoppingCGLS}
\end{align}
where $\tilde{\mathbf{x}}^{k}$ is the solution of \eqref{KroneckerRepresentationReducedLeastSquares}-\eqref{normalSystem} computed at the $k$-th iteration of the CGLS, and $\tilde{\mathbf{x}}^{0}$ is an initial guess that we chose as a zero vector. We stop the CGLS method when $\eta$ is below $10^{-6}$. Once the approximate solution has been computed with a prescribed tolerance, we quantify its (relative) accuracy by:
\begin{align}
\varepsilon=\left\|\tilde{\underline{X}}-\underline{X}_{T} \right\|_{2}/\left\| \underline{X}_{T} \right\|_{2}
\label{approximateSolutionAccuracy}
\end{align}
where $ \tilde{\underline{X}}$ is an approximate solution. Figure \ref{fig:Error}(a) shows the dependence of the accuracy \eqref{approximateSolutionAccuracy} on the bandwidth of the approximate solution.  The results are generated for $N=200$. As the bandwidth increases, the accuracy improves, as expected. It has been observed that the number of the iterations of the CGLS increases as the bandwidth increases (number of iterations to reach the stopping criteria defined by $\eta$). For example, for the bandwidth of $20$ it is $45$, and for the bandwidth of $300$ it is $235$. The condition number of $\underline{A}$ is $39$.  Figure \ref{fig:Error}(b) shows the row $600$ of the true and the approximate solution, and the error between them, computed for the bandwidth of $200$.
Next, we illustrate the dependence of the accuracy \eqref{approximateSolutionAccuracy} on the parameter $q$ in \eqref{method1LyapunovApproximateFinal}. We generate the results for $N=250$, bandwidth of $140$ in the iteration \eqref{sparsifiedChebyshev}, and the Chebyshev order of $M=20$. The approximation defined in \eqref{method1LyapunovApproximateFinal} has a total bandwidth of $294$. The results are shown in Fig. \ref{fig:ErrorConvergence}(a). As expected, the accuracy is improved by increasing $q$ and it confirms the exponential dependence predicted by \eqref{exponentialDecreaseFirstMethod}. In Fig. \ref{fig:ErrorConvergence}(b) we show the accuracy dependence on the number of iterations of the gradient projection method \eqref{steepestDescentMethod}. The initial guess is computed on the basis of  \eqref{method1LyapunovApproximateFinal} for $q=30$.
Finally, we test the computational and memory complexities of the proposed methods. We vary $N$, and measure the time necessary to compute the approximate solutions. The results are generated for the bandwidth of $150$. We also compare the complexity of the proposed methods with the complexity of the MATLAB function $\text{lyap}\left(\cdot \right)$. In order to compare the two methods, we chose their parameters such that their accuracy \eqref{approximateSolutionAccuracy} is approximately equal and below $0.03$ (for larger $N$ it takes a lot of time or even it is impossible to compute the "true" solution, so these parameters are determined heuristically, such that the accuracy is guaranteed for $N\le 600$). In the case of the first method, we stop the CGLS iteration when $\eta < 10^{-6}$. In the case of the second method, the parameters are: $q=60$, $M=20$ and the gradient projection iteration \eqref{steepestDescentMethod} is stopped after 50 iterations. Figure \ref{fig:Complexity}(a) shows the computational complexities of the first and second methods, together with the computational complexity of the function $\text{lyap}\left(\cdot\right)$. Figure \ref{fig:Complexity}(b) shows the memory complexity. 
From \ref{fig:Complexity}(a) it can be observed that the first method has the lowest computational complexity. The second method has a higher computational complexity than the first method, however, its memory complexity is lower. Computational and memory complexities of the first and the second method scale approximately linearly with $N$. This enables us to compute the approximate solution for larger values of $N$. For figure clarity, we presented results for $N$ up to $1500$. However, we tested the second method for problems up to $N=10^{4}$. The second method needs less than $50$ minutes to compute the solution for $N=10^{4}$. Due to the fact that we are not able to compute the "true" solution for such a large problem, we are not able to precisely quantify the accuracy of this approximate solution, but we expect that the accuracy is bounded by \eqref{exponentialDecreaseFirstMethod}. Low memory complexity of the second method allows us to compute the solution even for much larger $N$ (by extrapolating the results we estimate that for $N=10^{5}$, the second method can compute the approximate solution in less than $8$ hours). Finally, from Figure \ref{fig:Complexity} we can observe the $O(N^3)$ computational and $O\left(N^{2}\right)$ memory complexities of the function $\text{lyap}\left(\cdot\right)$. This function has the highest computational and memory complexities.
\subsection{Second example: randomly generated $\underline{A}$}
\begin{figure}[H]
  \centering
 \includegraphics[scale=0.25,trim=0mm 0mm 0mm 0mm ,clip=true]{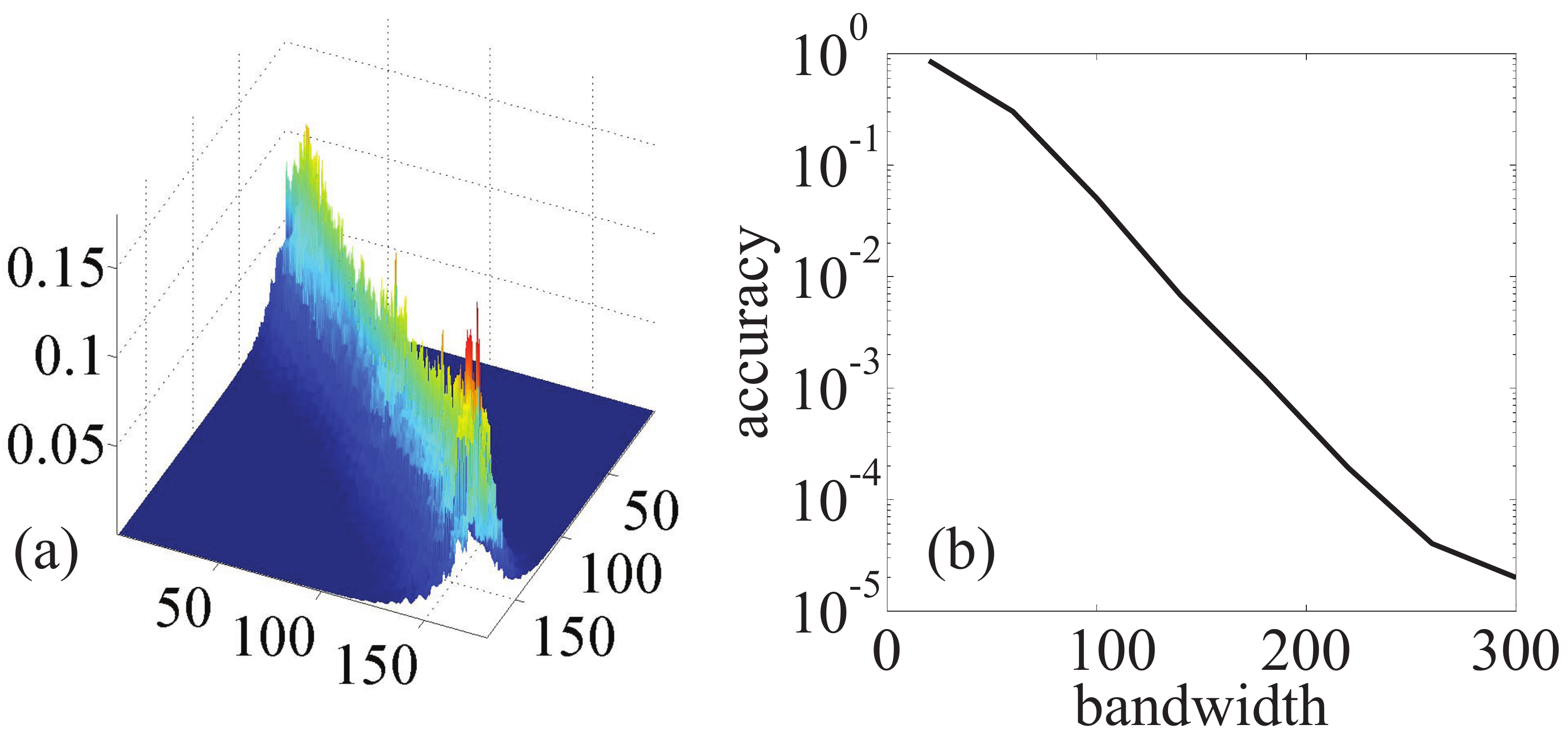}
 \caption{Randomly generated model. (a) Surface plot of $\underline{X}_{T}$, $N=30$. (b) Accuracy of the first method, $N=200$.}
\label{fig:SurfErrorRandom}
\end{figure}

Next, we illustrate the accuracy of the first method on a randomly generated model. Using the MATLAB function $\text{rand}\left(\cdot \right)$, we generate the matrices $A_{i,j}$, $j=i-1,i,i+1$ as a $6\times 6$ random matrices. After constructing the block tri-diagonal matrix $\underline{A}$ from these matrices, we define the following matrix $\underline{W}=\frac{1}{2}\left(\underline{A}+\underline{A}^{T}\right)+\nu I$, where the parameter $\nu$ has been chosen such that the matrix $\underline{W}$ is asymptotically stable. The matrix $\underline{P}$ is defined in the previous example. The surface plot of $\underline{X}_{T}$ for the pair of the coefficient matrices $(\underline{W},\underline{P})$ and for $N=30$, is shown in Fig. \ref{fig:SurfErrorRandom}(a).
From Fig. \ref{fig:SurfErrorRandom}(a) we see that the off-diagonal decay rate of $\underline{X}$ is fast, which directly follows from the fact that the matrix $\underline{W}$ is well-conditioned (its condition number is approximately $52$). Finally, in Fig. \ref{fig:SurfErrorRandom}(b) we show the accuracy dependence on the bandwidth, that improves as the bandwidth increases. The second method gives similar results and for the sake of brevity is omitted.
\subsection{Third example: 3D heat equation}
\label{thirdExpampleSection}
\begin{figure}[H]
  \centering
 \includegraphics[scale=0.22,trim=0mm 0mm 0mm 0mm ,clip=true]{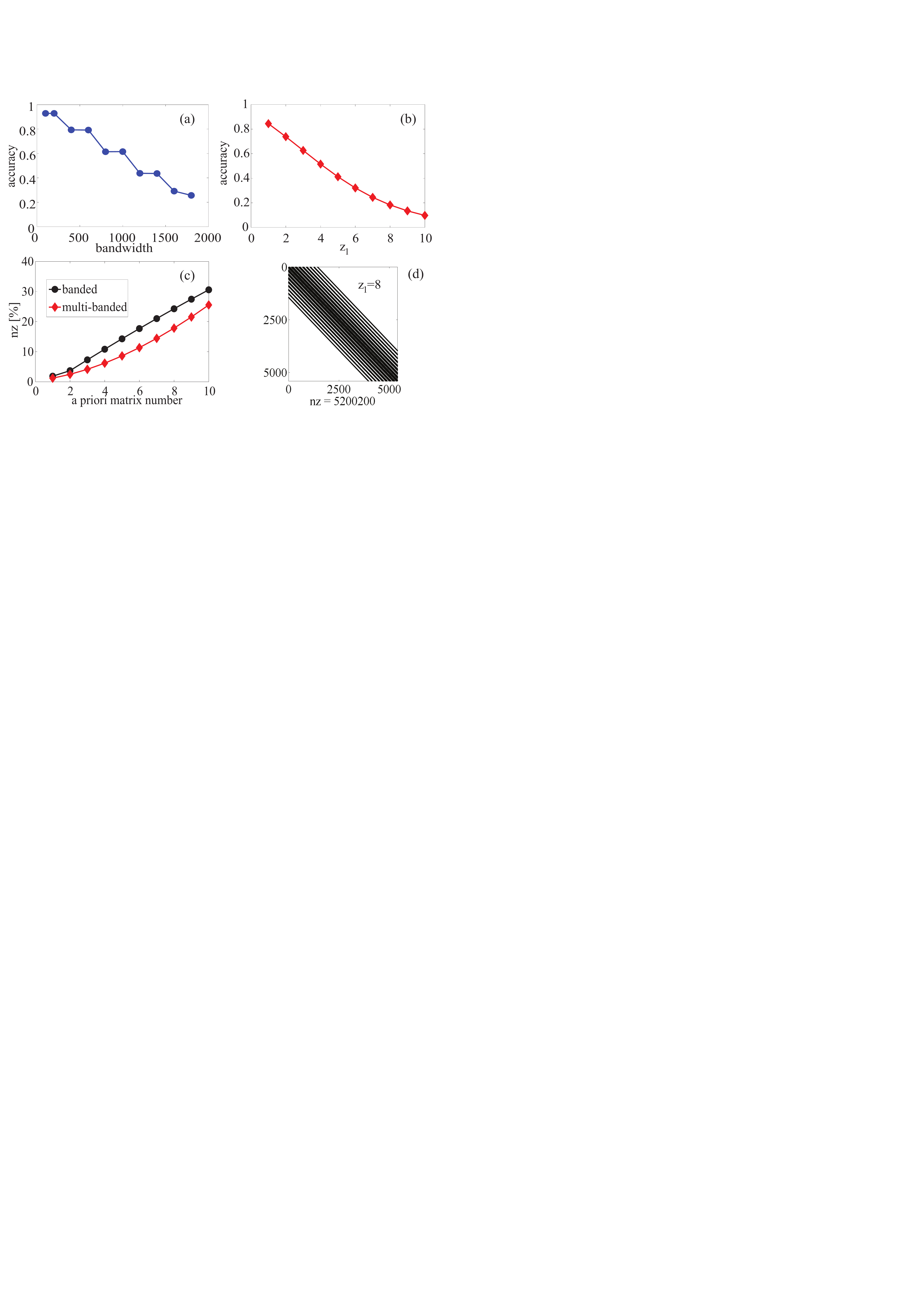}
 \caption{Discretized 3D heat equation. (a) Accuracy \eqref{approximateSolutionAccuracy} as a function of bandwidth for strictly banded a priori pattern. (b) Accuracy for the a-priori pattern computed using the method of Section \ref{generalSparsityPattern} (multi-banded pattern shown in (d) for $z_{1}=9$ ). (c) Comparison of the percentage of the non-zero elements (denoted by "nz") of a-priori sparsity patterns used to generate plots (a) and (b). (d) Plot of the multi-banded a-priori sparsity pattern, whose performance is shown in (b), the pattern is generated using the method in Section \ref{generalSparsityPattern} for $z_{1}=8$.}
\label{fig:AccErr3Dheat}
\end{figure}

We consider a 3D heat equation describing temperature change of a rectangular glass plate (diffusivity constant $3.4\times 10^{-7}$) used in optical systems \cite{haberThesis}. The heat equation is discretized using the finite difference method with the spatial discretization step of $0.001$ (all the units are in the SI system). The dimension of the grid in the $z$ direction is fixed to $6$, whereas in the $x-y$ plane(s) it is defined by an $N_{1}\times N_{1}$, where $N_{1}$ is the parameter that is varied. In total the grid has $6\times N_{1}^{2}$ points.  Such a discretization grid can be seen as an interconnection of $N=N^{2}_{1}$ subsystems $S_{i,j}$, where each subsystem's state consists of the temperatures in the $z$-direction: $\{T_{i,j,1},T_{i,j,2},\ldots,T_{i,j,6}\}$. That is, the local order of each $S_{i,j}$ is $n=6$, for more details see Chapter 2 of \cite{haberThesis}. 
\begin{figure}[H]
  \centering
 \includegraphics[scale=0.25,trim=0mm 0mm 0mm 0mm ,clip=true]{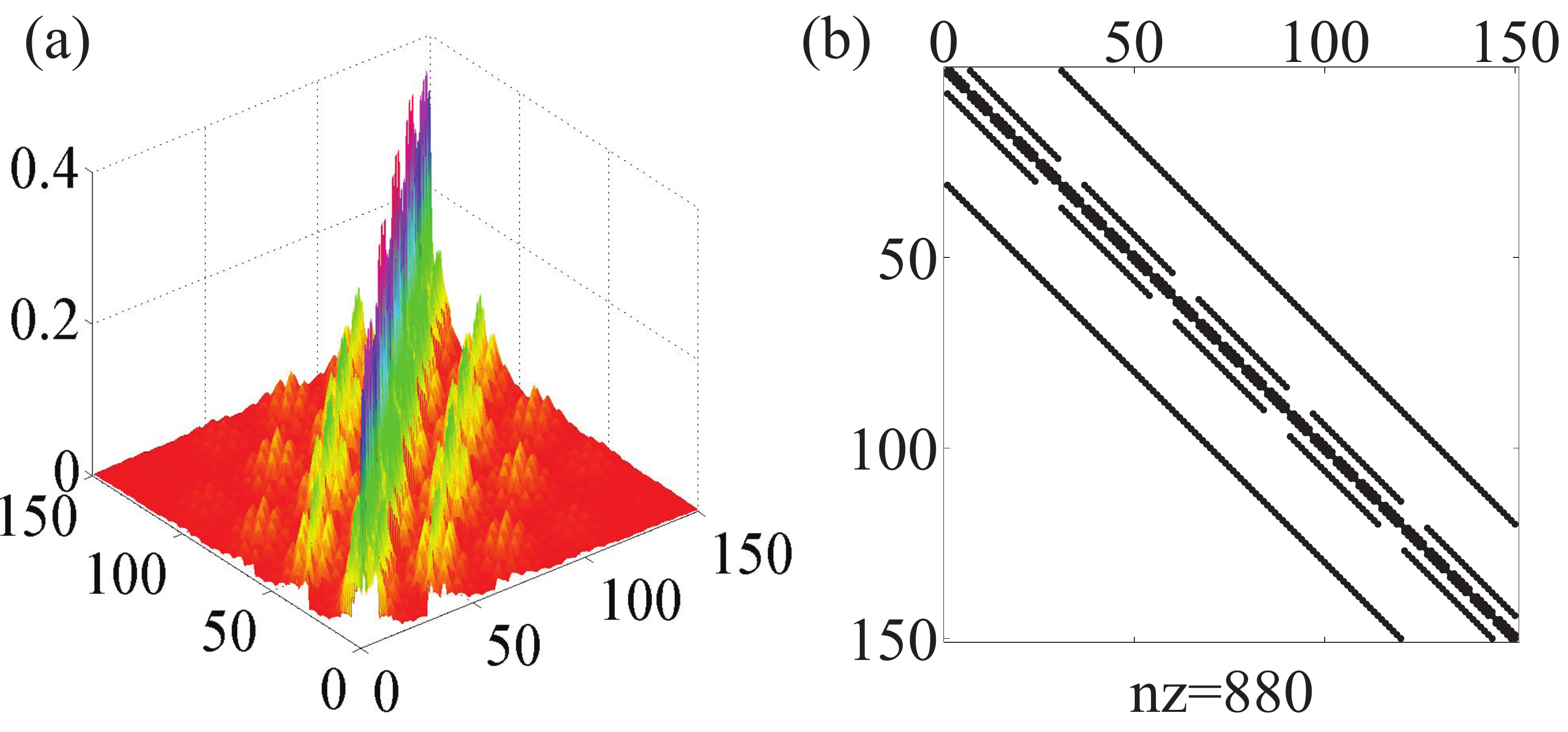}
 \caption{Discretized 3D heat equation. (a) Surface plot of $\underline{X}_{T}$. (b) Sparsity pattern of $\underline{A}$. $N=25$.}
\label{fig:SurfSparsity3Dheat}
\end{figure}
The matrix $\underline{P}$ is equal to the matrix used in the first example. The surface plot of $\underline{X}_{T}$ and sparsity pattern of $\underline{A}$ are shown in Fig. \ref{fig:SurfSparsity3Dheat}(a) and \ref{fig:SurfSparsity3Dheat}(b), respectively, for $N=25$, ($N_{1}=5$). Similarly to the first example, from Fig. \ref{fig:SurfSparsity3Dheat}(a) we see that $\underline{X}_{T}$ shows an off-diagonally decaying, oscillatory behavior. 
\\
We quantify the accuracy \eqref{approximateSolutionAccuracy} of the first method (results obtained by the second method are similar). For the CGLS we use the tolerance $\eta=10^{-6}$. We construct the model for $N_{1}=30$ which gives in total $N=900$ local subsystems ($\underline{A}\in \mathbb{R}^{5400\times 5400}$). The condition number of $\underline{A}$ is $72$ and it increases as $N$ is increased. Figure \ref{fig:AccErr3Dheat}(a) shows the accuracy for an a priori pattern equal to a banded matrix. As expected, the accuracy improves as the bandwidth increases. However, the improvement is slower compared to the first example, see Fig. \ref{fig:Error}(a). This is because $\underline{X}_{T}$ for the third example, has more dominant off-diagonal peaks than in the first example, as it can be observed in Fig. \ref{fig:SurfSparsity3Dheat}(a). Similarly to the first example, we noticed that the number of CGLS iterations increases (iterations to reach the prescribed tolerance) as the bandwidth is increased. Figure \ref{fig:AccErr3Dheat}(b) shows the accuracy for an a priori patten equal to a multi-banded matrix computed using the methodology proposed in Section \ref{generalSparsityPattern}. Such an a priori pattern is shown in Fig. \ref{fig:AccErr3Dheat}(d) for $z_{1}=8$. Finally, Fig. \ref{fig:AccErr3Dheat}(c) shows the percentage of the non-zero elements of the approximate solution for the two types of patterns. From Figures \ref{fig:AccErr3Dheat}(a)-(c) we conclude that multi-banded a priori pattern achieves better accuracy with a smaller number of non-zero elements compared to the purely banded a priori pattern. 
\section{Conclusions}
\label{sectionConclusions}
We considered the problem of computing a banded, approximate solution of the Lyapunov equation \eqref{LyapunovEquation}, in which the coefficient matrices are banded, symmetric matrices. We analyzed how the condition number of the coefficient matrix $\underline{A}$ influences the decay rate of the solution $\underline{X}$. We showed that the decay rate is faster if the condition number of $\underline{A}$ is smaller. We proposed two computationally efficient methods for approximating the solution by a banded matrix. Our results indicate that for a well-conditioned, sparse banded matrix $\underline{A}$, it is possible to compute a sparse, banded approximate solution of the Lyapunov equation with $O(N)$ complexity.
In the future work, the proposed methods will be used to develop computationally efficient methods for approximating the solution of the Riccati equation by a sparse matrix. 

\section*{Appendix}

\subsection{Computational complexity analysis of the second method}
\label{complexitySecond}
Due to the fact that we deal with (sparse) banded matrices, to estimate the computational complexity, we will not count multiplications and additions/subtractions of an arbitrary number with zero, because depending on the implementation (such as MATLAB sparse matrix toolbox), such operations do not require computing power. Let us assume that we have two banded, $Nn\times Nn$ matrices $X_{1}$ and $X_{2}$ with bandwidths $d_{1}$ and $d_{2}$, respectively. The matrix resulting from their multiplication $X_{3}=X_{1}X_{2}$ is a banded matrix with the bandwidth of $d_{3}=d_{1}+d_{2}$. The rows of $X_{3}$: $d_{3}/2+2,\ldots, Nn-d_{3}/2-1$ have $d_{1}+d_{2}+1$ non-zero entries, and remaining rows have smaller number of entries. For simplicity, we will assume that every row of $X_{3}$ has $d_{1}+d_{2}+1$ non-zero entries. To compute every entry on such a row, it takes maximum $d_{m}+1$ multiplications and $d_{m}$ additions, where $d_{m}=\min(d_{1},d_{2})$.
So in total, the number of operations that are necessary to compute $X_{3}$ is smaller than $Nn(d_{1}+d_{2}+1)(2d_{m}+1)$ operations. Under the same simplifications, we conclude that the number of operations necessary to compute $X_{4}=X_{1}\pm X_{2}$  is smaller than $Nn(d_{m}+1)$, and number of operations necessary to multiply a scalar $c$ and the matrix $X_{1}$ is smaller or equal to $Nn(d_{1}+1)$.
\\
Let us now turn our attention to the Chebyshev approximation of the matrix exponential.
The complexity of computing $\underline{A}_{1}$ using \eqref{matrixFunctions1} is much smaller than the complexity of subsequent operations and for simplicity it will be ignored. 
Consider the iteration \eqref{sparsifiedChebyshev} in which $T_{k}$ has a bandwidth of $d$ (produced by applying the operator $\mathcal{D}\left( \cdot \right)$ in the iteration $k-1$) and the matrix $\underline{A}$ has a bandwidth of $m$, $d>m$. It takes less than $Nn(m+1)+Nn(d+m+1)(2m+1)$ operations to compute $2\underline{A}_{1}T_{k}$ (multiplication of $\underline{A}_{1}$ by a scalar and multiplication of two banded matrices) plus additional $Nn(d+1)$ operations to compute $2\underline{A}_{1}T_{k}-T_{k-1}$. That is, the number of operations to compute one iteration \eqref{sparsifiedChebyshev} is less than $o_{1}=Nn(m+d+2+(d+m+1)(2m+1)   )$. For simplicity, we will ignore the operation count of the operator $\mathcal{D}\left( \cdot \right)$, because it depends on the implementation details. Due to the fact that we need $M-2$ iterations to compute all the matrices $T_{k}$ (we ignore complexity of computing $T_{1}$ and $T_{2}$), and because we need to multiply them with a scalar and add them together, the total number of operations to compute the approximation \eqref{truncatedChebyshevSeries}, is smaller than $o_{2}=(M-2)o_{1}+2NnM(d+1)$. The resulting matrix has the bandwidth of $d$.
\\
Let us now focus on \eqref{method1LyapunovApproximateFinal}. Assuming that the bandwidth $l$ of $\underline{P}$ is smaller than $d$, and using similar reasoning it can be concluded that the number of iterations is smaller than
\begin{align}
o_{3}=(2q+1)(o_{2}+Nn(d+l+1)(2l+1) \notag \\
+Nn(2d+l+1)(2d+1)+Nn(d+1)) \notag 
\end{align}
The important conclusion from the above analysis is that the complexity approximately scales linearly with $Nn$, $M$ and $q$, while on the other hand it scales quadratically with $d$.
\\
Similar analysis can be performed for the gradient projection method. However, it is hard to give a relatively precise estimate, because the number of step-size evaluations \eqref{ArmijoRule} in each iteration, cannot be predicted a priori. However, assuming that the number of step-size evaluations is very small, and using analysis similar to the Chebyshev approximation, it can be concluded that the computational complexity of every step of the gradient projection method scales linearly with $Nn$ and $d_{1}$ (there is no need to perform multiplications of two banded matrices of the bandwidths of $d_{1}$).

\bibliographystyle{unsrt}
\bibliography{bibl}

\begin{thebibliography}{10}

\bibitem{haber2014Gramian}
A.~Haber and M.~Verhaegen.
\newblock {Sparse Approximate Inverses of Gramians and Impulse Response
  Matrices of Large-scale Interconnected Systems}.
\newblock {\em {arXiv preprint arXiv:1405.2580v1 [cs.SY]}}, 2014.

\bibitem{pakazad2014distributed}
S.~K. Pakazad, A.~Hansson, M.~S. Andersen, and A.~Rantzer.
\newblock {Distributed Robustness Analysis of Interconnected Uncertain Systems
  Using Chordal Decomposition}.
\newblock {\em arXiv preprint arXiv:1402.2066}, 2014.

\bibitem{dandrea2003}
R.~D'Andrea and G.~Dullerud.
\newblock {Distributed Control Design for Spatially Interconnected Systems}.
\newblock {\em IEEE Transactions on Automatic Control}, 48:1478--1495, 2003.

\bibitem{jovbamTAC05platoons}
M.~R. Jovanovi\'c and B.~Bamieh.
\newblock {On the Ill-posedness of Certain Vehicular Platoon Control Problems}.
\newblock {\em IEEE Trans. Automat. Control}, 50(9):1307--1321, September 2005.

\bibitem{bamjovmitpat12}
B.~Bamieh, M.~R. Jovanovi\'c, P.~Mitra, and S.~Patterson.
\newblock {Coherence in Large-Scale networks: Dimension Dependent Limitations
  of Local Feedback}.
\newblock {\em IEEE Trans. Automat. Control}, 57(9):2235--2249, September 2012.

\bibitem{bamieh02distributedcontrol}
B.~Bamieh, O.~Paganini, and M.~A. Dahleh.
\newblock {Distributed Control of Spatially Invariant Systems}.
\newblock {\em IEEE Transactions on Automatic Control}, 47:1091--1118, 2002.

\bibitem{khan2008distributing}
U.~A. Khan and J.~M.~F. Moura.
\newblock {Distributing the Kalman Filter for Large-Scale Systems}.
\newblock {\em Signal Processing, IEEE Transactions on}, 56(10):4919--4935,
  2008.

\bibitem{gorinevsky2008design}
D.~Gorinevsky, S.~Boyd, and G.~Stein.
\newblock {Design of Low-Bandwidth Spatially Distributed Feedback}.
\newblock {\em Automatic Control, IEEE Transactions on}, 53(1):257--272, 2008.

\bibitem{motee2008}
N.~Motee and A.~Jadbabaie.
\newblock {Optimal control of spatially distributed systems}.
\newblock {\em IEEE Transactions on Automatic Control}, 53:1616--1629, 2008.

\bibitem{motee2013measuring}
N.~Motee and Q.~Sun.
\newblock {Measuring Sparsity in Spatially Interconnected Systems}.
\newblock In {\em Decision and Control (CDC), 2013 IEEE 52nd Annual Conference
  on}, pages 1520--1525. IEEE, 2013.

\bibitem{siami2014graph}
M.~Siami and N.~Motee.
\newblock {Graph-Theoretic Bounds on Disturbance Propagation in Interconnected
  Linear Dynamical Networks}.
\newblock {\em arXiv preprint arXiv:1403.1494}, 2014.

\bibitem{matni2014communication}
N.~Matni.
\newblock {Communication Delay Co-Design in $\mathcal{H}_{2}$ Distributed
  Control Using Atomic Norm Minimization}.
\newblock {\em arXiv preprint arXiv:1404.4911}, 2014.

\bibitem{andersen2014}
M.S. Andersen, S.K. Pakazad, A.~Hansson, and A.~Rantzer.
\newblock {Robust Stability Analysis of Sparsely Interconnected Uncertain
  Systems}.
\newblock {\em Automatic Control, IEEE Transactions on}, 59(8):2151--2156, Aug
  2014.

\bibitem{zhou2014}
T.~Zhou.
\newblock {Coordinated One-Step Optimal Distributed State Prediction for a
  Networked Dynamical System}.
\newblock {\em Automatic Control, IEEE Transactions on}, 58(11):2756--2771, Nov
  2013.

\bibitem{zhou2015controllability}
T.~Zhou.
\newblock {On the Controllability and Observability of Networked Dynamic
  Systems}.
\newblock {\em Automatica}, 52:63--75, 2015.

\bibitem{benner2004}
P.~Benner.
\newblock {Solving Large-Scale Control Problems}.
\newblock {\em Control Systems, IEEE}, 24(1):44 -- 59, 2004.

\bibitem{haberThesis}
A.~Haber.
\newblock {\em {Estimation and Control of Large-Scale Systems With an
  Application to Adaptive Optics for EUV Lithography}}.
\newblock PhD thesis, Delft University of Technology, Delft, The Netherlands,
  2014.

\bibitem{motter2013}
J.~Sun and A.~E. Motter.
\newblock {Controllability Transition and Nonlocality in Network Control}.
\newblock {\em Phys. Rev. Lett.}, 110:208701, May 2013.

\bibitem{gajic2008lyapunov}
Z.~Gajic and M.~T.~J. Qureshi.
\newblock {\em {Lyapunov Matrix Equation in System Stability and Control}}.
\newblock Courier Corporation, 2008.

\bibitem{mehrmann1991autonomous}
V.~L. Mehrmann.
\newblock {\em {The Autonomous Linear Quadratic Control Problem: Theory and
  Numerical Solution}}, volume 163.
\newblock Springer, 1991.

\bibitem{simoncini2013computational}
V.~Simoncini.
\newblock {Computational Methods for Linear Matrix Equations}.
\newblock {\em Technical Report. Dipartimento di Matematica. Universita di
  Bologna.}, 2013.

\bibitem{bini2012numerical}
D.~A. Bini, B.~Iannazzo, and B.~Meini.
\newblock {\em {Numerical Solution of Algebraic Riccati Equations}}, volume~9.
\newblock SIAM, 2012.

\bibitem{feitzinger2009inexact}
F.~Feitzinger, T.~Hylla, and E.~W. Sachs.
\newblock {Inexact Kleinman-Newton method for Riccati equations}.
\newblock {\em {SIAM Journal on Matrix Analysis and Applications}},
  31(2):272--288, 2009.

\bibitem{wang2014inexact}
X.~Wang, W.-W. Li, and L.~Dai.
\newblock {On Inexact Newton Methods Based on Doubling Iteration Scheme for
  Symmetric Algebraic Riccati Equations}.
\newblock {\em Journal of Computational and Applied Mathematics}, 260:364--374,
  2014.

\bibitem{benner2008numerical}
P.~Benner, J.~R. Li, and T.~Penzl.
\newblock {Numerical Solution of Large-Scale Lyapunov Equations, Riccati
  Equations, and Linear-Quadratic Optimal Control Problems}.
\newblock {\em Numerical Linear Algebra Appl}, 15(9):755--777, 2008.

\bibitem{bini2008fast}
D.~A. Bini, B.~Iannazzo, and F.~Poloni.
\newblock {A Fast Newton's Method for a Nonsymmetric Algebraic Riccati
  Equation}.
\newblock {\em {SIAM Journal on Matrix Analysis and Applications}},
  30(1):276--290, 2008.

\bibitem{benner2013numerical}
P.~Benner and J.~Saak.
\newblock {Numerical Solution of Large and Sparse Continuous Time Algebraic
  Matrix Riccati and Lyapunov Equations: A State of the Art Survey}.
\newblock Technical report, 2013.

\bibitem{Benner2011}
P.~Benner and H.~Fa{\ss}bender.
\newblock {On the Numerical Solution of Large-Scale Sparse Discrete-Time
  Riccati Equations}.
\newblock {\em {Adv. Comput. Math.}}, 35(2-4):119--147, 2011.

\bibitem{linfarjovTAC13admm}
F.~Lin, M.~Fardad, and M.~R. Jovanovi\'c.
\newblock {Design of Optimal Sparse Feedback Gains Via the Alternating
  Direction Method of Multipliers}.
\newblock {\em IEEE Trans. Automat. Control}, 58(9):2426--2431, September 2013.

\bibitem{linfarjovTAC11al}
F.~Lin, M.~Fardad, and M.~R. Jovanovi\'c.
\newblock {Augmented {L}agrangian Approach to Design of Structured Optimal
  State Feedback Gains}.
\newblock {\em IEEE Trans. Automat. Control}, 56(12):2923--2929, December 2011.

\bibitem{schuler2011}
S.~Schuler, P.~Li, J.~Lam, and F.~Allg\"{o}wer.
\newblock {Design of Structured Dynamic Output-Feedback Controllers for
  Interconnected Systems}.
\newblock {\em International Journal of Control}, 84(12):2081--2091, 2011.

\bibitem{benzi2007}
M.~Benzi and N.~Razouk.
\newblock {Decay Bounds and O(n) Algorithms for Approximating Functions of
  Sparse Matrices}.
\newblock {\em Electronic Transactions on Numerical Analysis}, pages 16--39,
  2007.

\bibitem{Haber:14subspace}
A.~Haber and M.~Verhaegen.
\newblock {Subspace Identification of Large-Scale Interconnected Systems}.
\newblock {\em Automatic Control, IEEE Transactions on}, 59(10):2754--2759,
  2014.

\bibitem{demko1984}
S.~Demko, W.~F. Moss, and P.~W. Smith.
\newblock {Decay Rates for Inverses of Band Matrices}.
\newblock {\em Mathematics of Computation}, 43(168):491--499, 1984.

\bibitem{benzi2015decay}
M.~Benzi and V.~Simoncini.
\newblock {Decay Bounds for Functions of Matrices with Banded or Kronecker
  Structure}.
\newblock {\em arXiv preprint arXiv:1501.07376}, 2015.

\bibitem{canuto2014decay}
C.~Canuto, V.~Simoncini, and M.~Verani.
\newblock {On the Decay of the Inverse of Matrices that are Sum of Kronecker
  Products}.
\newblock {\em Linear Algebra and its Applications}, 452:21--39, 2014.

\bibitem{simoncini2015lyapunov}
V.~Simoncini.
\newblock {The Lyapunov Matrix Equation. Matrix Analysis from a Computational
  Perspective}.
\newblock {\em arXiv preprint arXiv:1501.07564}, 2015.

\bibitem{verhaegen2007}
M.~Verhaegen and V.~Verdult.
\newblock {\em Filtering and System Identification: A Least Squares Approach}.
\newblock Cambridge University Press, 2007.

\bibitem{laub2005matrix}
A.~J. Laub.
\newblock {\em {Matrix Analysis for Scientists and Engineers}}.
\newblock Siam, 2005.

\bibitem{bertsekas1999nonlinear}
D.~P. Bertsekas.
\newblock {\em {Nonlinear programming}}.
\newblock Athena Scientific, 1999.

\bibitem{shao2014finite}
M.~Shao.
\newblock {On the Finite Section Method for Computing Exponentials of
  Doubly-Infinite Skew-Hermitian Matrices}.
\newblock {\em Linear Algebra and its Applications}, 451:65--96, 2014.

\bibitem{benzi1999bounds}
M.~Benzi and G.~H. Golub.
\newblock {Bounds for the Entries of Matrix Functions with Applications to
  Preconditioning}.
\newblock {\em {BIT Numerical Mathematics}}, 39(3):417--438, 1999.

\bibitem{huckle1999apriori}
T.~Huckle.
\newblock {Approximate Sparsity Patterns for the Inverse of a Matrix and
  Preconditioning}.
\newblock {\em Applied numerical mathematics}, (30):291--303, 1999.

\bibitem{pan1991improved}
V.~Pan and R.~Schreiber.
\newblock {An Improved Newton Iteration for the Generalized Inverse of a Matrix
  with Applications}.
\newblock {\em SIAM Journal on Scientific and Statistical Computing},
  12(5):1109--1130, 1991.

\bibitem{saad2003iterative}
Y.~Saad.
\newblock {\em Iterative Methods for Sparse Linear Systems}.
\newblock Society for Industrial and Applied Mathematics, 2003.

\bibitem{bjorck1996numerical}
A.~Bj{\"o}rck.
\newblock {\em Numerical methods for least squares problems}.
\newblock Siam, 1996.

\bibitem{chow1998approximate}
E.~Chow and Y.~Saad.
\newblock {Approximate Inverse Preconditioners via Sparse-Sparse Iterations}.
\newblock {\em SIAM Journal on Scientific Computing}, 19(3):995--1023, 1998.

\bibitem{chow2000priori}
E.~Chow.
\newblock {A Priori Sparsity Patterns for Parallel Sparse Approximate Inverse
  Preconditioners}.
\newblock {\em SIAM Journal on Scientific Computing}, 21(5):1804--1822, 2000.

\bibitem{benzi1998sparse}
M.~Benzi and M.~Tuma.
\newblock {A Sparse Approximate Inverse Preconditioner for Nonsymmetric Linear
  Systems}.
\newblock {\em SIAM Journal on Scientific Computing}, 19(3):968--994, 1998.

\bibitem{bergamaschi2000efficient}
L.~Bergamaschi and M.~Vianello.
\newblock {Efficient Computation of the Exponential Operator for Large, Sparse,
  Symmetric Matrices}.
\newblock {\em Numerical linear algebra with applications}, 7(1):27--45, 2000.

\bibitem{mason2003chebyshev}
J.~C. Mason and D.~C. Handscomb.
\newblock {\em {Chebyshev Polynomials}}.
\newblock Chapman \& Hall/CRC, 2003.

\bibitem{Mathar2006}
R.~J. Mathar.
\newblock {Chebyshev Series Expansion of Inverse Polynomials}.
\newblock {\em J. Comput. Appl. Math.}, 196(2):596--607, November 2006.

\bibitem{bergamaschi2003efficient}
L.~Bergamaschi, M.~Caliari, and M.~Vianello.
\newblock {Efficient Approximation of the Exponential Operator for Discrete 2D
  Advection--Diffusion Problems}.
\newblock {\em Numerical linear algebra with applications}, 10(3):271--289,
  2003.

\bibitem{grasedyck2003solution}
L.~Grasedyck, W.~Hackbusch, and B.N. Khoromskij.
\newblock {Solution of Large Scale Algebraic Matrix Riccati Equations by Use of
  Hierarchical Matrices}.
\newblock {\em Computing}, 70(2):121--165, 2003.

\bibitem{li2012weighted}
Z.-Y. Li and Y.~Wang.
\newblock {Weighted Steepest Descent Method for Solving Matrix Equations}.
\newblock {\em International Journal of Computer Mathematics},
  89(8):1017--1038, 2012.

\bibitem{benner2013efficient}
P.~Benner, J.~Saak, M.~Stoll, and H.~K. Weichelt.
\newblock {Efficient Solution of Large-Scale Saddle Point Systems Arising in
  Riccati-based Boundary Feedback Stabilization of Incompressible Stokes Flow}.
\newblock {\em SIAM Journal on Scientific Computing}, 35(5):S150--S170, 2013.

\bibitem{lehoucq1998arpack}
R.~B. Lehoucq, D.~C. Sorensen, and C.~Yang.
\newblock {\em {ARPACK Users' Guide: Solution of Large-scale Eigenvalue
  Problems with Implicitly Restarted Arnoldi Methods}}, volume~6.
\newblock SIAM, 1998.

\bibitem{haber2013OL}
A.~Haber, A.~Polo, S.~K. Ravensbergen, H.~P. Urbach, and M.~Verhaegen.
\newblock {Identification of a Dynamical Model of a Thermally Actuated
  Deformable Mirror}.
\newblock {\em Opt. Lett.}, 38(16):3061--3064, Aug 2013.

\bibitem{haber2013predictive}
A.~Haber, A.~Polo, I.~Maj, S.F. Pereira, H.P. Urbach, and M.~Verhaegen.
\newblock {Predictive Control of Thermally Induced Wavefront Aberrations}.
\newblock {\em Optics express}, 21(18):21530--21541, 2013.

\bibitem{Haber:13}
A.~Haber, A.~Polo, C.~S. Smith, S.~F. Pereira, P.~Urbach, and M.~Verhaegen.
\newblock {Iterative Learning Control of a Membrane Deformable Mirror for
  Optimal Wavefront Correction}.
\newblock {\em Appl. Opt.}, 52(11):2363--2373, Apr 2013.

\end{thebibliography}
\end{document}